\RequirePackage{fix-cm}
\documentclass[smallextended]{svjour3}       %
\smartqed  %
\usepackage{graphicx}
\DeclareGraphicsExtensions{.png,.pdf,.jpg}
\usepackage[colorlinks=true]{hyperref}

\usepackage[cmex10]{amsmath}
\usepackage{amssymb}
\usepackage{multirow}
\usepackage{tabularx}

\usepackage[T1]{fontenc}
\usepackage{newtxtext} %

\usepackage{graphics} %
\usepackage{epsfig} %
\graphicspath{{images/}}
\usepackage{graphicx}
\usepackage{mwe}
\usepackage{url}
\usepackage{comment}
\usepackage[cmex10]{amsmath} %
\usepackage{booktabs}
\usepackage{amssymb}  %
\usepackage{xcolor}

\usepackage{algorithm}

\usepackage{algpseudocode}

\usepackage{accents}

\newcommand{\BesselI}[2]{\mathcal{I}_{#1} \left( #2 \right)}

\newcommand{\EXP}[1]{\exp{\left( #1 \right)}}
\newcommand{\EXPSHORT}[1]{\mathrm{e}^{ #1 }}

\newcommand{\BEGMATRIX}[1]{\left[\begin{array}{#1}}
	\newcommand{\ENDMATRIX}{\end{array}\right]}

\newcommand{\COS}[1]{\cos{\left( #1 \right)}}
\newcommand{\SIN}[1]{\sin{\left( #1 \right)}}

\def\det{\mathop{\rm det}}

\newcommand{\ba}{\mathbf{a}}
\newcommand{\bb}{\mathbf{b}}
\newcommand{\bc}{\mathbf{c}}
\newcommand{\bd}{\mathbf{d}}
\newcommand{\be}{\mathbf{e}}

\newcommand{\bp}{\mathbf{p}}

\newcommand{\bs}{\mathbf{s}}
\newcommand{\bt}{\mathbf{t}}
\newcommand{\bu}{\mathbf{u}}
\newcommand{\bv}{\mathbf{v}}

\newcommand{\bw}{\mathbf{w}}

\newcommand{\bz}{\mathbf{z}}

\newcommand{\bB}{\mathbf{B}}

\newcommand{\bF}{\mathbf{F}}

\newcommand{\bH}{\mathbf{H}}
\newcommand{\bI}{\mathbf{I}}

\newcommand{\bL}{\mathbf{L}}
\newcommand{\bM}{\mathbf{M}}

\newcommand{\bR}{\mathbf{R}}

\newcommand{\bU}{\mathbf{U}}
\newcommand{\bV}{\mathbf{V}}

\newcommand{\bY}{\mathbf{Y}}

\newcommand{\bmu}{\boldsymbol{\mu}}

\newcommand{\bomega}{\boldsymbol{\omega}}

\newcommand{\bLambda}{\boldsymbol{\Lambda}}
\newcommand{\bPhi}{\boldsymbol{\Phi}}

\newcommand{\bOmega}{\boldsymbol{\Omega}}

\newcommand{\boldzero}{\boldsymbol{0}}

\newcommand{\transpose}{^\mathrm{T}}
\newcommand{\hermitian}{^\mathrm{H}}

\newcommand{\littletaller}{\mathchoice{\vphantom{\big|}}{}{}{}}
\newcommand\restr[2]{{%
  \left.\kern-\nulldelimiterspace %
  #1 %
  \littletaller %
  \right|_{#2} %
  }}

\newcommand{\newcommenter}[2]{%
  \expandafter\newcommand\csname #1\endcsname[1]{%
    \fcolorbox{#2}{#2}{\color{white}\textsf{\textbf{#1}}}
    {\color{#2}##1}}%
  \expandafter\newcommand\csname at#1\endcsname{%
    \fcolorbox{#2}{#2}{\color{white}\textsf{\textbf{@#1}}}
    {\color{#2}}}%
  \expandafter\newcommand\csname #1cite\endcsname[1]{%
    \csname #1\endcsname {[##1]}
  }%
  \expandafter\newcommand\csname #1ref\endcsname[1]{%
    \csname #1\endcsname {$\blacktriangleright$##1}
  }%
  \expandafter\newcommand\csname #1hl\endcsname[2]{%
    \colorbox{#2}{\color{white}\textsf{\textbf{#1}}}\sethlcolor{Azure2}\hl{##2}~%
    \expandafter\ifx\csname commentarrow\endcsname\relax$\leftarrow$\else \commentarrow[#2]\fi~%
    {\color{#2}##1}}%
  \expandafter\newcommand\csname #1st\endcsname[2]{%
    \colorbox{#2}{\color{white}\textsf{\textbf{#1}}}\sout{##2}~%
    \expandafter\ifx\csname commentarrow\endcsname\relax$\leftarrow$\else \commentarrow[#2]\fi~%
    {\color{#2}##1}}%
}

\makeatletter
\@ifpackageloaded{xcolor}{}{%
\usepackage[table,11names,dvipsnames,svgnames]{xcolor}%
}
\makeatother

\usepackage{colortbl}

\definecolor{darkgreen}{RGB}{0,100,0}
\definecolor{darkred}{RGB}{100,0,0}
\definecolor{darkblue}{RGB}{100,0,0}
\definecolor{darkcyan}{RGB}{100,100,0}

\newtheorem{thm}{Theorem}
\newtheorem{prop}[thm]{Proposition}
\usepackage{easybmat} %
\newcommand{\bmat}[1]{\begin{bmatrix}#1 \end{bmatrix}}

\usepackage{mathdots}

\newcommand{\bphat}{\hat{\bp}}

\newcommand{\ADJ}[2]{\mathrm{ad}_{#1}^{#2}}

\newcommenter{ernesto}{darkgreen}
\newcommenter{todo}{darkred}
\newcommenter{dario}{darkblue}
\newcommenter{francesco}{darkcyan}

\usepackage{svg}

\newcounter{mysubequations}
\renewcommand{\themysubequations}{(\roman{mysubequations})}
\newcommand{\mysubnumber}{\refstepcounter{mysubequations}\themysubequations}
\makeatletter\def\makeheadbox{}\makeatother

\begin{document}

\title{A Lie-Algebraic Formulation of Real Spherical Harmonic Rotation%
}

\author{Ernesto Fontana         \and
        Dario Lodi Rizzini      \and
        Francesco Pasetti
}

\institute{E. Fontana, D. Lodi Rizzini \at
              Robotics and Intelligent Machines Laboratory (Rimlab), Dipartimento di Ingegneria e Architettura, and Centro Interdipartimentale per l'Energia e l'Ambiente (CIDEA), University of Parma, 43124 Parma, Italy \\
              \email{ernesto.fontana@unipr.it, dario.lodirizzini@unipr.it}
           \and
           F. Pasetti \at
              Robotics and Intelligent Machines Laboratory (Rimlab), Dipartimento di Ingegneria e Architettura, University of Parma, 43124 Parma, Italy \\
              \email{francesco.pasetti@unipr.it}
}
\date{Preprint submitted to the \emph{Journal of Scientific Computing}}

\maketitle

\begin{abstract}
This work presents a Lie-algebraic formulation of the rotation of real spherical harmonic (RSH) coefficients, together with an algorithm for computing the corresponding rotation matrix. 
From the differential operators realizing infinitesimal rotations on the sphere, we derive the generators of the $\mathfrak{so}(3)$ representation acting on the degree-$\ell$ RSH subspace, and obtain the RSH-rot matrix through the exponential map. 
As an alternative that avoids a general-purpose matrix exponential, we factorize RSH-rot into anti-diagonal form, so that the rotation angle enters only through an elementary rotation about the $z$ axis, while an orthogonal similarity matrix depends only on the rotation axis.
The similarity matrix is characterized by two independent routes: through the adjoint action on the Lie algebra, and from the eigenvectors of the RSH-rot matrix, which are computed in real arithmetic by a re-orthogonalized ladder recurrence.
The formulation is free of the singularities affecting Euler-angle methods and clarifies the algebraic structure of RSH-rot.
As the Lie algebra corresponds to the tangent space at the identity, RSH-rot can be linearized, making rotation-dependent functions of RSH series explicitly differentiable, with gradients and Hessians in closed form.
Numerical experiments against the Ivanic--Ruedenberg recursion show that the Lie-algebraic methods remain accurate to round-off at high degrees, where the recursion loses accuracy, that single degrees are computed faster, that many rotations about a common axis are evaluated at about the cost of one, and that the closed-form derivatives are exact and cheaper than finite differences.
\keywords{Real Spherical Harmonics \and Spherical Harmonic Rotation \and Lie algebra $\mathfrak{so}(3)$ \and Matrix exponential \and Rotation estimation}
\subclass{33C55 \and 22E70 \and 65F60 \and 65F15}
\end{abstract}

\section{Introduction}
\label{sec:introduction}

Spherical harmonics (SH) form an orthogonal basis for representing functions on the sphere and play a fundamental role across applied mathematics and science, with applications ranging from the solution of partial differential equations~\cite{DassiosFokas08}, to quantum mechanics~\cite{VarshalovichMoskalevKhersonskii88}, physical chemistry~\cite{RicoLopezEmaRamirez04}, computer graphics~\cite{Sloan08}, 3D shape analysis and perception~\cite{MakadiaDaniilidis06,AlthloothiMahoorVoyles13}, neural networks~\cite{PassaroZitnick23}, communication~\cite{AlemKhalidKennedy16}.
Under proper conditions, any well-defined function $f$ on the unit sphere can be expanded into a series of orthogonal SH, which provide a complete orthogonal basis for square-integrable functions on the sphere.
Their orthogonality enables efficient projection and computation in the spectral domain, while their eigenfunction structure provides a natural decomposition into angular modes and simplifies differential, convolution, and rotational operations.

SH are usually given in two forms: \emph{complex-valued SH} (CSH) and \emph{real-valued SH} (RSH).
While CSH are commonly adopted in physics and mathematical formulations, RSH are often preferred when representing real functions on the unit sphere, particularly in applications such as chemistry and molecular modeling~\cite{Nikolaev11} and computer graphics~\cite{FridovichKeilYuTancikChenRechtKanazawa22}.
Under a rotation $\bR$, a spherical harmonic of degree $\ell$ is transformed into a linear combination of the $2 \ell + 1$ spherical harmonics of the same degree.
Consequently, the action of $\bR$ on the degree-$\ell$ SH subspace can be represented by a $(2 \ell + 1) \times (2 \ell + 1)$ rotation matrix $\bM_\ell(\bR)$, whose elements depend only on the rotation.
We call \emph{RSH-rot} such representation of rotation on degree-$\ell$ coefficients of RSH. 

The computation of RSH-rot is required in several operations on the unit sphere such as the computation of integrals~\cite{BourhraraCollignonAylajVidal26} or interpolation~\cite{DellAccioTommaso25}.
The expression of the rotation matrix for the complex case was originally provided by Wigner~\cite[p. 158--161]{Wigner59} through the so-called D-matrix. 
Gimbutas and Greengard~\cite{GimbutasGreengard09} rotate spherical harmonic expansions through a factorization into shear and diagonal operators, which is fast and stable but confined to CSH and yielding no rotation matrix or Lie-algebraic structure to differentiate.
Other works have provided algorithms for the rotation matrix of RSH, such as those of Blanco et al.~\cite{Blanco97} or Pinchon et al.~\cite{Pinchon07}. 
Both construct SH-rot from the Euler angle representation of the rotation, so their evaluation requires three successive angle-dependent factors. 
Such parameterizations inherit the drawbacks of Euler angles, most notably the gimbal lock singularity, and the accumulated trigonometric factors also limit accuracy at high degrees $\ell$.
These methods derive $\bM_\ell(\bR)$ from the Euler angle representation of the 3D rotation $\bR$, which suffers from inconveniences such as gimbal lock.
The method proposed by Ivanic and Ruedenberg~\cite{IvanicRuedenberg96} recursively computes the RSH rotation from the lower-order matrix.
Other works derive the RSH-rot for some low-order $\ell$ from the rotation described in angle-axis form~\cite{Romanowski07,Maintz16}.
Some works~\cite{KrivanekKonttinenPattanaikBouatouchZara06,BigiFrauxBrowningCeriotti23} propose fast approximations of RSH-rot for practical applications. 

In this paper, we present a novel formulation of RSH-rot grounded in the Lie algebra $\mathfrak{so}(3)$ and derive from it an algorithm for computing the rotation matrix $\bM_\ell(\bR)$.
We begin by deriving the generators of the representation of $\mathfrak{so}(3)$ acting on the degree-$\ell$ RSH subspace directly from the differential operators that realize infinitesimal rotations on the sphere.
The Lie algebra of RSH-rot is obtained as a linear combination of the generators, from which the RSH-rot matrix follows by the exponential map.
As an alternative to a general-purpose matrix exponential, we derive an algorithm based on the similarity between the given Lie algebra element and the anti-diagonal generator $\bL_{\ell z}$ associated to rotations around the $z$ axis.
The RSH-rot is then computed as the product of this similarity matrix and the RSH-rot $\bM_{\ell,z}$ about the $z$ axis, so that the rotation angle enters only through an elementary rotation and the similarity matrix depends only on the axis.
We characterize the similarity matrix by two independent routes.
The first is based on the action of the rotation on the Lie algebra itself through the adjoint representation.
The second computes the eigenvectors of the RSH-rot matrix through a ladder recurrence in real arithmetic.
We prove that the similarity matrix depends directly on the real and imaginary parts of the eigenvectors, and we show that the recurrence must be re-orthogonalized to be numerically stable.

The proposed approach offers several advantages.
First, it is free of the singularities that affect methods built on Euler angles, and it computes each degree independently of the lower ones.
Second, it gives a deeper understanding of the algebraic structure of RSH-rot.
Third, the separation of axis and angle in the anti-diagonal factorization makes many rotations about the same axis inexpensive, and turns the correlation of two RSH expansions along that axis into a trigonometric polynomial.
Fourth, the Lie algebra permits the linearization of RSH-rot through the tangent space, so that rotation-dependent functions of RSH series have gradients and Hessians in closed form.
Numerical experiments compare the proposed methods with the recursion of Ivanic and Ruedenberg in terms of computation time and accuracy.
They show that the Lie-algebraic methods remain accurate to round-off at high degrees, where the recursion loses accuracy, while the recursion remains the fastest method when all degrees of a single rotation are required.

In summary, the contributions of this paper are the following.
\begin{itemize}
	\item A Lie-algebraic characterization of RSH-rot, yielding the generators of the $\mathfrak{so}(3)$ representation from differential operators and the RSH-rot matrix itself through the exponential map.
	\item An anti-diagonal factorization of RSH-rot that separates the rotation axis from the angle, with the similarity matrix obtained both from the action on the Lie algebra and from the eigenvectors of the RSH-rot matrix and its Lie-algebra representation, computed by a stable real recurrence.
	\item Explicit differentiability of rotation-dependent functions of RSH series, enabling the use of RSH-rot in gradient-based optimization.
\end{itemize}

\section{RSH Rotation Matrix}
\label{sec:rsh_rot_definition}

\subsection{Definition of RSH Rotation}

The main subject of this paper are RSH and their rotation. 
According to the standard convention~\cite{Bayin18} normalized RSH are special functions defined as 
\begin{equation}
  Y_{\ell m}(\theta, \varphi)
  {=}\sqrt{2-\delta_{m0}}
  \sqrt{\frac{2 \ell{+}1}{4 \pi}\frac{(\ell{-}|m|)!}{(\ell{+}|m|)!}}
    P_\ell^{|m|}(\cos\theta)\Phi_m(\varphi)
  \label{eqn:rsh_definition}
\end{equation}
where $\theta$ and $\varphi$ are respectively the colatitude and longitude in the unit sphere, $\delta_{\cdot \cdot}$ is the Kronecker delta, $P_\ell^m(\cos\theta)$ is the associated Legendre function without the Condon--Shortley phase, and $\Phi_m(\varphi)$ is the sinusoidal function
\begin{equation}
  \Phi_m(\varphi) = \left\{\begin{array}{ll}
   \SIN{|m| \varphi} & -\ell \leqslant m < 0 \\ 
   1                     & m = 0 \\
   \COS{m \varphi}   & 0 < m \leqslant \ell 
   \end{array}\right. 
   \label{eqn:rsh_definition_sinusoid}
\end{equation}
The positive integer $\ell$ is called degree, order or band of the SH. 
The SH in band $\ell$ can be written as degree $\ell$ homogeneous polynomials in Cartesian coordinates.  
The secondary index $m$ is an integer varying $-\ell \leqslant m \leqslant \ell$. 

In this paper, we adopt two specific notations. 
First, the argument of SH is usually described with an element of the unit sphere $\bphat \in \mathbb{S}^2$ instead of using the colatitude and longitude angles $\theta$ and $\varphi$: $\bphat = [\sin\theta \cos\varphi, \sin\theta \sin\varphi, \cos\theta]\transpose$.
This notation is not only more compact, but it also suggests that the element of the domain are unit vectors that can be rotated. 
Second, the RSHs of the band $\ell$ are collected into the vector function $\bY_\ell(\bphat) = [Y_{\ell,-\ell}(\bphat), \ldots, Y_{\ell,\ell}(\bphat)]\transpose$. 
As discussed in the following, the action of rotation on SHs affects together the functions on the same band $\ell$. 
Such vector notation $\bY_\ell(\bphat)$ allows better identification of matrix.

One important application is the expansion of square integrable functions on the unit sphere $f \in L^2(\mathbb{S}^2)$ defined  into a series of SH:
\begin{equation}
  f(\bphat) = \sum_{\ell \geqslant 0} \bc_{\ell}\transpose \bY_{\ell}(\bphat)
  \label{eqn:func_rsh_series}
\end{equation}
where $\bc_{\ell} = [c_{\ell,-\ell}, \dots, c_{\ell,\ell}]\transpose \in \mathbb{R}^{2\ell+1}$ is the vector of real RSH coefficients for band $\ell$. 
This usage is equivalent to Fourier series expansion for functions defined on unit circle $\mathbb{S}^1$. 
The SH provides an orthonormal basis that enables operations among functions supporting such representation.  
The orthogonality of RSH is defined in relation to integral product of two real functions defined on the unit sphere $f, g \in L^2(\mathbb{S}^2)$ and is  
\begin{equation}
  \langle f(\cdot), g(\cdot) \rangle 
  = \int_{\mathbb{S}^2} f(\bphat) ~g(\bphat) d\bphat
\end{equation}
This operator $\langle\cdot,\cdot\rangle$ is also applied to vector and matrix valued functions in the sense that it is computed on each element of the vector or of the matrix.   

\begin{property}[Orthogonality of RSH]
    Two RSHs $Y_{\ell m}$ and $Y_{\ell' m'}$ are orthogonal with respect to the integral product, i.e.
    \begin{equation}
    \langle Y_{\ell m}(\cdot), Y_{\ell' m'}(\cdot) \rangle = \delta_{\ell \ell'} ~\delta_{m m'}.
    \end{equation}
\end{property}
Orthogonality allows to compute the series coefficients in eq.~\eqref{eqn:func_rsh_series} as $\bc_{\ell} = \langle f(\cdot), \bY_{\ell}(\cdot) \rangle$. 

This property can be exploited to evaluate how rotation acts on each RSH. 
Given a matrix $\bR \in SO(3)$, the function $f(\bphat) = Y_{\ell m}(\bR \bphat)$ can be expanded into RSH series as in eq.~\eqref{eqn:func_rsh_series}.
It can be proved that the coefficients of such series $\bc_{\ell'} = \boldzero$ when $\ell' \neq \ell$ and, therefore, $Y_{\ell m}(\bR \bphat)$ can be written as a linear combination of the spherical harmonics $Y_{\ell m'}(\bphat)$ with $m'=-\ell, \ldots, \ell$.  
In other words, SH of degree $\ell$ form an irreducible representation of the rotation group $SO(3)$~\cite{Selig10}, and therefore rotations preserve the $(2\ell+1)$-dimensional subspace spanned by the harmonics of that degree.
In summary, the rotated vector of RSH functions of order $\ell$ can be written as 
\begin{equation}
	\bY_\ell(\bR \bphat) = \bM_\ell(\bR) ~\bY_\ell(\bphat)
	\label{eqn:rsh_rot_linear}
\end{equation}
The matrix $\bM_\ell(\bR)$ in the above equation is the main subject of this paper. 
For compatibility with SH indices, the row and column indices of $\bM_\ell$ as well as of all the $(2\ell+1)$-square matrices appearing in this paper vary from $-\ell$ to $+\ell$ instead of the common positive indices. 
In this paper, we also adopt a minimal representation of rotation in the form of the so called \emph{Rodrigues vectors}. 
A Rodrigues vector is vector $\ba \in \mathbb{R}^3$ such that $\ba = \theta_a \hat{\ba}$ where $\hat{\ba}$ is the unit rotation axis $\|\hat{\ba}\| = 1$ and $\theta_a \geqslant 0$ is the rotation angle. 
This vector encodes the axis-angle representation of rotation and also Lie algebra parameters for Cartesian basis. 
RSH-rot matrices are sometimes referred to as $\bM_\ell(\bR_\ba)$ or $\bM_{\ell,\ba}$ using such compact parameters encoding axis-angle notation. 

\begin{definition}[RSH Rotation Matrix]
    The matrix $\bM_\ell(\bR)$ in eq.~\eqref{eqn:rsh_rot_linear} collecting the coefficients of the linear combination is called \emph{RSH rotation} (RSH-rot). 
\end{definition}

The orthogonality property enables rewriting the RSH-rot in a convenient form.
The integral product between $\bY_{\ell}(\bR \bu)$ and $\bY_\ell\transpose(\bu)$, after the substitution of eq.~\eqref{eqn:rsh_rot_linear}, yields
\begin{equation}
	\bM_\ell(\bR) 
	= \langle \bY_\ell(\bR \cdot), \bY_\ell\transpose(\cdot) \rangle
    = \int_{\mathbb{S}^2} \bY_\ell(\bR \bphat) ~\bY_\ell\transpose(\bphat) ~d\bphat
	\label{eqn:rsh_rot_def}
\end{equation}
Such expression does not enable direct computation of $\bM_\ell(\bR)$, but can be used to derive its properties.

\subsection{Properties of RSH Rotation}

This section illustrates relevant properties of the RSH-rot matrix directly derived from its definition. 

\begin{prop}
	\label{thm:rsh_rot_composition_eigen}
	Let $\bR_\ba, \bR_\bb, \bR_\bc \in SO(3)$ be three rotation matrices with Rodrigues parameters $\ba$, $\bb$ and $\bc$. 
	Then, the following relations hold
	\begin{equation}
		\begin{aligned}\setcounter{mysubequations}{0}
			\mysubnumber\quad & 
			\bM_\ell(\bI_3) = \bI_{2\ell+1} \\
			\mysubnumber\quad & 
			\bM_\ell(\bR_{\bb} \bR_{\ba}) = \bM_\ell(\bR_{\bb}) \bM_\ell(\bR_{\ba}) \\
            \mysubnumber\quad & 
            \bM_\ell(\bR_{\bc}) (\bM_\ell(\bR_{\bb}) \bM_\ell(\bR_{\ba})) = (\bM_\ell(\bR_{\bc}) \bM_\ell(\bR_{\bb})) \bM_\ell(\bR_{\ba}) \\
			\mysubnumber\quad & 
			\bM_\ell(\bR_{\ba}^{-1}){=}\bM_\ell(\bR_{\ba}\transpose) = \bM_\ell\transpose(\bR_{\ba}) = \bM_\ell^{-1}(\bR_{\ba}) 
		\end{aligned}
		\label{eqn:rsh_rot_composition_eigen}
	\end{equation}
\end{prop}

\begin{proof} 
	Property~$(i)$ follows from substitution of rotation $\bR=\bI_3$ in  \eqref{eqn:rsh_rot_def} and the orthogonality of SH. \\
	Property~$(ii)$ stems from observing that $\bY_{\ell}((\bR_{\bb} \bR_{\ba}) \bphat)$ has RSH-rot $\bM_{\ell}(\bR_{\bb} \bR_{\ba})$.
    Given $\bphat'=\bR_{\ba} \bphat$, applying eq.~\eqref{eqn:rsh_rot_linear} twice gives $\bY_{\ell}(\bR_{\bb} \bphat') = \bM_{\ell}(\bR_{\bb})\bY_{\ell}(\bphat') = \bM_{\ell}(\bR_{\bb})\bM_{\ell}(\bR_{\ba})\bY_{\ell}(\bphat)$. 
    Since the RSH of band $\ell$ are linearly independent, the matrix in eq.~\eqref{eqn:rsh_rot_linear} is unique, and the composition follows.\\
    Property~$(iii)$ is the direct result of associativity on $SO(3)$ rotation matrices and (ii).
	Property~$(iv)$ follows from the definition \eqref{eqn:rsh_rot_def} applied to $\bR\transpose$ and from a change of variable $\bv = \bR\transpose \bu$ as
	\begin{align}
		& \bM_\ell(\bR\transpose) = \int_{\mathbb{S}^2} \bY_{\ell}(\bR\transpose \bu) ~\bY_{\ell}\transpose(\bu) ~d\bu 
		= \int_{\mathbb{S}^2} \bY_{\ell}(\bv) ~\bY_{\ell}\transpose(\bR \bv) ~d\bv \nonumber\\
		&= \int_{\mathbb{S}^2} (\bY_{\ell}(\bR \bv) ~\bY_{\ell}\transpose(\bv))\transpose ~d\bv 
		= \left( \int_{\mathbb{S}^2} \bY_{\ell}(\bR \bv) ~\bY_{\ell}\transpose(\bv) ~d\bv \right)\transpose
		= \bM_\ell\transpose(\bR) \nonumber
	\end{align}
	The proof is completed by observing that $\bR_{\ba}^{-1}=\bR_{\ba}\transpose$ and by applying $(i)$ and $(ii)$ to $\bM_\ell(\bR_{\ba}^{-1} \bR_{\ba})$. 
\end{proof}  

Proposition~\ref{thm:rsh_rot_composition_eigen} demonstrates that RSH-rot matrices are representations of the Lie group of 3D rotations $SO(3)$. 
RSH-rot satisfy the axioms of associativity, identity, and inverse element. 
RSH-rot is also a map from a rotation to an action on vector space $\mathbb{R}^{2\ell+1}$. 
Like all the actions on vector space, it is represented in the form of a matrix of $GL(2\ell+1, \mathbb{R})$. 
In this case, $\bM_\ell(\bR)$ is orthogonal due to eq.~\eqref{eqn:rsh_rot_composition_eigen}(iv). 
Moreover, $\det\bM_\ell(\bR)=+1$, since $SO(3)$ is connected, $\bM_\ell(\cdot)$ is continuous, and $\bM_\ell(\bI_3)=\bI_{2\ell+1}$; hence $\bM_\ell(\bR)\in SO(2\ell+1)$. 
Moreover, if $\bM_\ell(\cdot)$ is a representation of a Lie group, there is also a corresponding representation $\bL_{\ell}(\cdot)$ of the Lie algebra. 
$\bL_{\ell}$ maps the elements of the algebra $so(3)$ to their corresponding $(2\ell+1)$-square matrices and preserves the Lie bracket operation in the new domain. 
Since $\bM_\ell \in SO(2\ell+1)$, $\bL_\ell \in so(2\ell+1)$ have the form of skew-symmetric matrices. 
The Lie algebra approach is the core contribution of this work and is discussed in section~\ref{sec:rsh_lie_algebra_generators}.  
Before that we complete the characterization of the RSH-rot. 

\begin{prop}
  \label{thm:rsh_band_l_rot_eigenvector}
  Let $\bR_\ba \in SO(3)$ be a rotation matrix with Rodrigues vector $\ba$, $\hat{\ba} \in \mathbb{S}^2$ be its unitary rotation axis, and $\bM_{\ell}(\bR_\ba)$ be the corresponding RSH-rot matrix for band $\ell  \geqslant 0$. 
  Then, $\bv = Y_{\ell }(\hat{\ba})$ is an eigenvector of $\bM_{\ell}(\bR_\ba)$ with eigenvalue $1$. 
\end{prop}
\begin{proof}
  Eq.~\eqref{eqn:rsh_rot_linear} writes the rotated RSH function $\bY_{\ell}(\bR \bphat)$ as the product of $\bM_{\ell}(\bR)$ and of $\bY_{\ell}(\bphat)$. 
  The axis $\hat{\ba}$ is an eigenvector of $\bR$, i.e. $\bR \hat{\ba} = \hat{\ba}$. 
  If we choose $\bphat = \hat{\ba}$, then  
  \begin{equation}
  	\bM_{\ell}(\bR) \bY_\ell(\hat{\ba})
  	= \bY_\ell(\bR \hat{\ba}) = \bY_\ell(\hat{\ba})
  \end{equation}
  Thus, $\bY_\ell(\hat{\ba})$ satisfies the definition of eigenvector with eigenvalue $1$ for matrix $\bM_\ell(\bR)$. 
\end{proof}

Proposition~\ref{thm:rsh_band_l_rot_eigenvector} relates the exact expression of the eigenvector associated to the only real eigenvalue. 
This expression is exploited for the computation of the other eigenvectors in section~\ref{sec:rsh_rot_eigenvectors}.
It also allows discussion about functions with axial symmetry and the coefficients of their RSH expansion.

\begin{prop}
    \label{thm:rsh_axis_symmetric}
	Let $f: \mathbb{S}^2 \to \mathbb{R}$ be a function expanded in RSH series with coefficients $\bc_{\ell}$ of each band $\ell \geqslant 0$ as in eq.~\eqref{eqn:func_rsh_series}. 
    Let $\bR_{\ba} \in SO(3)$ be a rotation associated with the Rodrigues vector $\ba = \theta_a \hat{\ba}$. 
    Then,  $f$ is symmetric with respect to axis $\hat{\ba}$, i.e. $f(\bR_{\theta_a \hat{\ba}} \bphat) = f(\bphat)$ for any possible $\theta_a$, iff 
	\begin{align}
		\bc_\ell &= k_\ell ~\bY_\ell(\hat{\ba})  \label{eqn:rsh_coefficients_symmetric}
	\end{align}
	for some constants $k_\ell \in \mathbb{R}$. 
	Moreover, the function $f$ can be written using Legendre polynomials as 
	\begin{equation}
		f(\bphat) 
		= \sum_{l \geqslant 0} \frac{2\ell+1}{4 \pi} k_\ell~P_l(\hat{\ba} \transpose \bphat)
        \label{eqn:rsh_axis_symmetric_legendre}
	\end{equation}
\end{prop}

\begin{proof}
  Let $\ba = \theta_a \hat{\ba}$ for any $\theta_a$. The RSH expansions of original and rotated $f$ are respectively 
  \begin{align}
  	f(\bphat) 
  	&= \sum_{\ell \geqslant 0} \bc_\ell\transpose ~\bY_\ell(\bphat) \\
  	f(\bR_\ba \bphat) 
  	&= \sum_{\ell \geqslant 0} \bc_\ell\transpose (\bM_\ell(\bR_\ba) \;\bY_\ell(\bphat))  
  	 = \sum_{\ell \geqslant 0} (\bM_\ell\transpose(\bR_\ba) \bc_\ell)\transpose \;\bY_\ell(\bphat)
  \end{align}
  \emph{Only if}. If $f$ is $\hat{\ba}$-symmetric and $\bR_\ba$ is a rotation around axis $\hat{\ba}$, then $f(\bphat) = f(\bR_\ba \bphat)$. 
  The coefficients of the two functions must be equal for every $\ell $ and $m$ and 
  \begin{align}
  	\bc_\ell &= (\bM_\ell (\bR_\ba))\transpose \bc_\ell = \bM_\ell (\bR_\ba\transpose) \bc_\ell
  \end{align}
  Thus, $\bc_\ell$ is an eigenvector of $\bM_\ell(\bR_\ba\transpose)$ for eigenvalue $1$. 
  Choosing an angle $\theta_a$ for which the eigenvalue $1$ is simple (Proposition~\ref{thm:rsh_eigenvalues}), Proposition~\ref{thm:rsh_band_l_rot_eigenvector} gives $\bc_\ell = k_\ell\bY_\ell(\hat{\ba})$ ($\hat{\ba}$ is the rotation axis of both $\bR_\ba$ and its inverse $\bR_\ba\transpose$). 
  \\
  \emph{If}. If the RSH coefficients are those in eq.~(\ref{eqn:rsh_coefficients_symmetric}), then 
  \begin{align}
  	f(\bR_\ba \bphat) 
    &= \sum_{\ell \geqslant 0} (\bM_\ell\transpose(\bR_\ba) ~k_\ell \bY_\ell(\hat{\ba}))\transpose \;\bY_\ell(\bphat) \\
  	&= \sum_{\ell  \geqslant 0} (k_\ell \bY_\ell(\hat{\ba})\transpose)  \bY_\ell(\bphat) 
     = \sum_{\ell  \geqslant 0} \bc_\ell\transpose  \bY_\ell(\bphat) 
     = f(\bphat) 
  \end{align}
  Proposition~\ref{thm:rsh_band_l_rot_eigenvector} has been used to remove dependency from $\bM_\ell\transpose(\bR_\ba)$. 
  Finally, the addition theorem of SH~\cite[Eq.~(10.360)]{Bayin18} allows rewriting as in eq.~\eqref{eqn:rsh_axis_symmetric_legendre}. 
\end{proof}

The above result is useful for detecting symmetries from RSH coefficients. 

\subsection{RSH Rotation around Axis $\hat{\bz}$}
\label{sec:rsh_rot_z}

The main goal of this paper is the evaluation of RSH-rot and will be addressed in the general case.
Before concluding this section it is convenient to present its exact expression for the specific case of rotations around axis $\hat{\bz}$.  
Let $\bphat=[\sin\theta \cos\varphi, \sin\theta \sin\varphi, \cos\theta]\transpose$ be a point on unit sphere  identified by parameters $(\theta, \varphi)$ where $\theta$ is the colatitude and $\varphi$ the longitude. 
The rotation $\bR_{z}(\alpha)$ with axis $\hat{\bz}$ and angle $\alpha$ transforms $\bphat$ into $\bphat'$ with parameters $(\theta, \varphi{+}\alpha)$. 
When $\bR_{z}(\alpha)$ is applied to RSH, it affects only the term $\bPhi_m(\varphi)$ in eq.~\eqref{eqn:rsh_definition_sinusoid}. 
In particular, the rotated components $\Phi_{-m}(\varphi+\alpha)$ and $\Phi_{m}(\varphi+\alpha)$ (assume $m>0$) are equal to
\begin{equation}
  \begin{bmatrix} \Phi_{-m}(\varphi+\alpha) \\ \Phi_{m}(\varphi+\alpha) \end{bmatrix}
  =
  \begin{bmatrix} 
    \COS{m \alpha} & \SIN{m \alpha} \\ 
   -\SIN{m \alpha} & \COS{m \alpha} \
  \end{bmatrix}
  \begin{bmatrix} \Phi_{-m}(\varphi) \\ \Phi_{m}(\varphi) \end{bmatrix}
  \label{eqn:rhs_rot_z_phi_shift}
\end{equation}
The above equations enable writing the rotated RSH $Y_{\ell,-m}(\bR_{z}(\alpha) \bphat)$ and $Y_{\ell,m}(\bR_{z}(\alpha) \bphat)$. 
In order to orderly write the complete RSH-rot $\bM_{\ell z}(\alpha) \triangleq \bM_{\ell}(\bR_{z}(\alpha))$ we adopt the following notation. 
Given a square matrix $\bM$ with dimension $2\ell+1$, the symbol $\bM^{[m]}$, $m \geqslant 0$, refers to the square sub-matrix with the $(2m+1) \times (2m+1)$ elements in the center of $\bM$, i.e. the elements with row and column indices in $\{-m, \ldots, 0, \ldots, +m\}$. 
The recurrent expression of $\bM_{\ell z}(\alpha)$ and its minors stems directly from eq.~\eqref{eqn:rhs_rot_z_phi_shift} with $m=1, \ldots, \ell$ as

\begin{equation}
  \begin{gathered}
\bM_{\ell z}^{[1]}(\alpha) =
\begin{bmatrix}
  \cos\alpha & 0 & \sin\alpha \\
  0          & 1 & 0 \\
 -\sin\alpha & 0 & \cos\alpha 
\end{bmatrix}
\\
\bM_{\ell z}^{[m]}(\alpha) =
\left[
\begin{array}{ccccc}
  \COS{m \alpha}     &        &        &        &  \SIN{m \alpha} \\[2pt]
       & \multicolumn{3}{c}{\multirow{3}{*}{%
           \fbox{\rule[-6pt]{0pt}{18pt}$\bM_{\ell z}^{[m{-}1]}( \alpha )$}}} &   \\[2pt]
       &        &        &        &        \\[2pt]
       &        &        &        &        \\[2pt]
 - \SIN{m \alpha}  &        &        &        &  \COS{m \alpha}
\end{array}
\right]
  \end{gathered}
  \label{eqn:rsh_rot_z}
\end{equation}

Eq.~\eqref{eqn:rsh_rot_z} enables the computation of the eigenvalues $\bM_{\ell z}(\alpha)$. 
The characteristic polynomial of the minor $\bM_{\ell z}^{[m]}(\alpha)$ is
\begin{align}
&|\lambda \bI_{2m{+}1}{-}\bM_{\ell z}^{[m]}(\alpha)| =
\left|
\left[
\begin{array}{ccccc}
  \COS{m \alpha}     &        &        &        &  \SIN{m \alpha} \\[2pt]
       & \multicolumn{3}{c}{\multirow{3}{*}{%
           \fbox{\rule[-6pt]{0pt}{18pt}$\lambda\bI_{2m-1}-M_{\ell z}^{[m{-}1]}( \alpha )$}}} &   \\[2pt]
       &        &        &        &        \\[2pt]
       &        &        &        &        \\[2pt]
 - \SIN{m \alpha}  &        &        &        &  \COS{m \alpha}
\end{array}
\right]
\right|\nonumber\\
&= (\lambda^2 - 2 \COS{m \alpha} \lambda + 1) ~|\lambda\bI_{2m-1}-M_{\ell z}^{[m{-}1]}( \alpha )| \nonumber\\
&= (\lambda - \EXPSHORT{i m \alpha}) (\lambda - \EXPSHORT{-i m \alpha}) 
   ~|\lambda\bI_{2m-1}-M_{\ell z}^{[m{-}1]}( \alpha )|
\end{align}
Thus, the complete characteristic polynomial of $\bM_{\ell z}(\alpha)$ is
\begin{equation}
  |\lambda \bI_{2\ell{+}1}{-}\bM_{\ell z}(\alpha)| 
  = (\lambda-1) \prod_{m=1}^{\ell} (\lambda - \EXPSHORT{i m \alpha}) (\lambda - \EXPSHORT{-i m \alpha})
\end{equation}
The eigenvalues of the matrix are $\EXPSHORT{-i \ell \alpha}, \ldots, 1, \ldots, \EXPSHORT{+i \ell \alpha}$. 
Any rotation matrix $\bR_{\ba} \in SO(3)$ with axis $\hat{\ba}$ and angle $\theta_a$ can be factorized into the product of $\bR=\bU\;\bR_{z}\;(\theta_a) \bU\transpose$ where $\bU$ is a proper orthonormal matrix that rotates $\hat{\ba}$ to $\hat{\bz}$. 
From Proposition~\ref{thm:rsh_rot_composition_eigen}(ii) follows that $\bM_\ell(\bR_{\ba})$ is similar to $\bM_{\ell z}(\theta_a)$ and, therefore, it has the same eigenvalues. 
This is stated in the following proposition. 

\begin{prop}
    \label{thm:rsh_eigenvalues}
    Let $\bR_{\ba} \in SO(3)$ be a rotation with axis $\hat{\ba}$ and angle $\theta_a$ and $\bM_\ell(\bR_{\ba})$ be its RSH-rot representation for band $\ell$.
    Then, the eigenvalues of $\bM_\ell(\bR_{\ba})$ are
	\begin{equation}
	    \EXPSHORT{-i \ell \theta_a}, \EXPSHORT{-i (\ell-1) \theta_a}, \ldots, 1, \ldots, \EXPSHORT{+i (\ell-1) \theta_a}, \EXPSHORT{+i \ell \theta_a}
	\end{equation}
\end{prop}

The factorization of $\bM_\ell$ and the characterization of its eigenvalues and eigenvectors is addressed in section~\ref{sec:rsh_rot_eigenvectors}, but first we present the analysis of its Lie algebra.

\section{Lie Algebra of RSH Rotation}
\label{sec:rsh_lie_algebra}

The RSH rotation matrices $\bM_\ell$ provide representations of the rotation group $SO(3)$ and are themselves orthogonal matrices in $SO(2\ell+1)$, as shown in Proposition~\ref{thm:rsh_rot_composition_eigen}. As for any Lie group, there is a corresponding representation $\bL_\ell$ of its Lie algebra (LA), $\mathfrak{so}(3)$. In particular,
$\bL_\ell: \mathfrak{so}(3) \to \mathfrak{so}(2\ell+1)$
is a Lie algebra homomorphism into the space of $(2\ell+1)\times(2\ell+1)$ skew-symmetric matrices. In particular, it preserves the Lie bracket: for any $\bOmega_\ba, \bOmega_\bb \in \mathfrak{so}(3)$,
\begin{equation}
  \bL_{\ell}([\bOmega_\ba, \bOmega_\bb])
  =
  [\bL_{\ell}(\bOmega_\ba), \bL_{\ell}(\bOmega_\bb)].
  \label{eqn:lie_bracket_representation}
\end{equation}
The representation $\bL_\ell$ provides a complementary description of the rotation group and offers several advantages.

First, the elements $\bOmega_\ba \in \mathfrak{so}(3)$ are represented by compact $3\times 3$ skew-symmetric matrices. Each such matrix can be expressed as a linear combination of three basis elements, which we refer to as \emph{generators}. The coefficients of this linear combination are collected in a vector $\ba\in\mathbb{R}^3$; the subscript in $\bOmega_\ba$ makes this relationship explicit. By linearity, the representation $\bL_\ell(\bOmega_\ba)$ is given by the same linear combination of the representations of the generators, with the same coefficients $\ba$. This observation provides a straightforward way to compute $\bL_\ell(\bOmega_\ba)$. 
Moreover, the exponential map relates the Lie algebra representation to the corresponding RSH rotation matrix:
$\bM_\ell = \EXPSHORT{\bL_\ell(\bOmega_\ba)}$.

Second, the Lie algebra describes the tangent space of the group at the identity and therefore provides a natural framework for differentiating RSH rotations. 
Several problems require estimating the rotation that optimizes an objective function defined over RSH representations, for example, the alignment of two directional distributions. 
In these cases, computing $\bM_\ell(\bR)$ directly from the Euler angles (or Euler parameters) of $\bR$ is not sufficient; a representation that supports differentiation with respect to the rotation parameters is also needed.

This section presents three main contributions. 
First, we derive a set of generators for the Lie algebra representation, which enables the computation of $\bM_\ell(\bR)$. 
Second, we establish several properties of the Lie bracket that are useful for the subsequent derivations.
Third, we present a factorization of the RSH-rot that is convenient for its computation.  

\subsection{Generators of the Lie Algebra of RSH Rotations}
\label{sec:rsh_lie_algebra_generators}

The Lie algebra $\mathfrak{so}(3)$ is a three-dimensional vector space. Consequently, each of its elements can be expressed as a linear combination of the generators associated with a chosen basis. A common choice is the basis whose elements correspond to rotations about the axes of the reference frame. In the defining, three-dimensional representation, these generators are the real skew-symmetric matrices
\begin{equation}
  \bOmega_x =
  \begin{bmatrix}
    0 & 0 &  0 \\
    0 & 0 & -1 \\
    0 & 1 &  0
  \end{bmatrix},
  \quad
  \bOmega_y =
  \begin{bmatrix}
    0 & 0 &  1 \\
    0 & 0 &  0 \\
   -1 & 0 &  0
  \end{bmatrix},
  \quad
  \bOmega_z =
  \begin{bmatrix}
    0 & -1 &  0 \\
    1 &  0 &  0 \\
    0 &  0 &  0
  \end{bmatrix},
  \label{eq:so3_cartesian_basis}
\end{equation}
which satisfy the commutation relations
$[\bOmega_k,\bOmega_j]=\varepsilon_{kjm}\bOmega_m$
for $k,j,m\in\{x,y,z\}$, where $\varepsilon_{kjm}$ is the Levi-Civita symbol. The set
$\mathcal{B}_{xyz}=\{\bOmega_x,\bOmega_y,\bOmega_z\}$
is referred to as the \emph{Cartesian basis}.

A second basis, which is more convenient when working with the complex extension of the algebra, is obtained through the change of basis
\begin{equation}
  \bOmega_+ = \bOmega_x + i\,\bOmega_y, \qquad
  \bOmega_- = \bOmega_x - i\,\bOmega_y, \qquad
  \bOmega_0 = \bOmega_z,
  \label{eq:so3_spherical_basis}
\end{equation}
and is referred to as the \emph{complex basis}, 
$\mathcal{B}_{\pm 0}=\{\bOmega_+,\bOmega_-,\bOmega_0\}$.
These generators satisfy $[\bOmega_k,\bOmega_j]=i \varepsilon_{kjm}\bOmega_m$.
In general, given a basis $\{\bOmega_k\}_k$ of the abstract algebra, any element $\bOmega_\ba\in\mathfrak{so}(3)$ can be written as
\begin{equation}
  \bOmega_\ba = \sum_{k=1}^{3} a_k\,\bOmega_k.
  \label{eq:so3_general_element}
\end{equation}
Here, the coefficient vector is $\ba=[a_1,a_2,a_3]\transpose$. When the Cartesian basis is used, $\ba$ has the geometric interpretation of a rotation vector (or Rodrigues vector), as introduced in Section~\ref{sec:rsh_rot_definition}. Its norm,
$\theta_a=\|\ba\|$, is the rotation angle, while the normalized vector
$\hat{\ba}=\ba/\theta_a$ specifies the rotation axis. Accordingly,
$\bOmega_\ba=\theta_a\bOmega_{\hat{\ba}}$.

The action of $\bOmega_\ba$ on the spherical harmonics (SH) of band $\ell$ is obtained by applying the representation $\bL_\ell(\cdot)$ term by term. Thus, the band-$\ell$ element associated with the rotation vector $\ba$ is
\begin{equation}
  \bL_{\ell,\ba}
  \triangleq \bL_\ell(\bOmega_\ba)
  = \sum_{k=1}^{3} a_k\,\bL_{\ell,k},
  \qquad
  \bL_{\ell,k}\triangleq\bL_\ell(\bOmega_k),
  \label{eqn:so3_band_ell_element}
\end{equation}
where $\bL_{\ell,k}$ denotes the representation of the $k$-th generator for band $\ell$. 
Thanks to the linearity, the Lie algebra representation can be scaled as $\bL_{\ell, \ba} = \theta_a \bL_{\ell, \hat{\ba}}$.
The Cartesian and complex bases of the represented algebra are, respectively,
$\mathcal{B}^{\ell}_{xyz}=\{\bL_{\ell x},\bL_{\ell y},\bL_{\ell z}\}$
and
$\mathcal{B}^{\ell}_{\pm 0}=\{\bL_{\ell +},\bL_{\ell -},\bL_{\ell 0}\}$.
Our ultimate goal is to derive the matrices associated with the Cartesian basis. To this end, however, it is convenient to work first with the complex basis.

The Lie algebra action, which can be interpreted intuitively as an infinitesimal rotation, can be described both by the skew-symmetric matrices introduced above and by differential operators. When applied to spherical harmonics, these differential operators describe the infinitesimal change of a function under rotation~\cite[Sec.~10.8]{Bayin18}. There is a direct correspondence between the differential operator $L_k$ acting on SH, the Lie algebra basis element $\bOmega_k$, and its band-$\ell$ representation $\bL_{\ell,k}$. The differential operators depend on the basis adopted. In terms of the angular variables $\theta$ and $\varphi$ used for complex spherical harmonics (CSH), they are defined as
\begin{equation}
  L_0 = -i \frac{\partial}{\partial\varphi},
  \qquad
  L_{\pm} = \EXPSHORT{\pm i\varphi}
  \left(
    \pm \frac{\partial}{\partial\theta}
    + i\cot\theta\frac{\partial}{\partial\varphi}
  \right).
  \label{eqn:operator_definition}
\end{equation}
Their action on CSH is given by~\cite[ Eqs.~(8.191)--(8.192)]{Bayin18}
\begin{equation}
  L_0 Y_{\ell}^{m} = mY_{\ell}^{m},
  \qquad
  L_+Y_{\ell}^{m} = \gamma_{\ell,m}Y_{\ell}^{m+1},
  \qquad
  L_-Y_{\ell}^{m} = \gamma_{\ell,-m}Y_{\ell}^{m-1},
  \label{eqn:operator_action_to_csh}
\end{equation}
where
$\gamma_{\ell,m}=\sqrt{(\ell-m)(\ell+m+1)}$.
The operator $L_+$ increases the index $m$ of the spherical harmonic on which it acts, whereas $L_-$ decreases it. 
For this reason, $L_+$ and $L_-$ are commonly referred to as \emph{ladder operators}. 
The increment or decrement occurs only when $m$ remains between the limits $-\ell$ and $\ell$ as we implicitly assume in the following.
Their coefficients satisfy the identities
$\gamma_{\ell,m}=\gamma_{\ell,-(m+1)}$
and
$\gamma_{\ell,-m}=\gamma_{\ell,m-1}$. 

We use the differential operators to derive the matrices associated with the Cartesian basis $\mathcal{B}^{\ell}_{xyz}$. Two aspects of the conversion from CSH to RSH must be considered. First, the definition in Eq.~\eqref{eqn:operator_definition} follows the angular-momentum convention, in which rotations are represented through the purely imaginary exponential $\EXPSHORT{-i\bL}$~\cite[Sec.~10.8]{Bayin18}. Since our objective is to obtain real matrix representations of the Lie algebra, an additional factor of $-i$ must therefore be taken into account. Second, the CSH $Y_\ell^m$ in Eq.~\eqref{eqn:operator_definition} must be replaced by the corresponding RSH $Y_{\ell,m}$, and the complex operators must be expressed in terms of the real operators $L_x$, $L_y$, and $L_z$. The conversion formulas between CSH and RSH are given in Section~\ref{sec:appendix_csh_conversion}. The complex and real operators are related by
$L_0=L_z$ and $L_\pm=L_x\pm iL_y$.

The action of $L_z=L_0$ on the RSH provides the matrix representation $\bL_{\ell z}$. For $m=0$ (omitting the argument $\bphat$ of SH) we obtain
\begin{equation}
   L_z Y_{\ell,0}=L_0Y_\ell^0=0.
   \label{eqn:diff_lz_rsh_0}
\end{equation}
It is convenient to discuss the complementary cases of positive and negative index $m$, when $m \neq 0$. 
Since the formula of conversion formulas between CSH and rSH are different, for clarity we assume that $0 < m \leqslant \ell$. 
The action of $L_{0}$ to $Y_{\ell, -m}$ and $Y_{\ell, m}$ has interlaced formulas that are handled together as  
\begin{align}
  L_{z} &\begin{bmatrix} Y_{\ell, -m} \\ Y_{\ell, m} \end{bmatrix}
  = \frac{1}{\sqrt2}
     \begin{bmatrix} 
        i & - i (-1)^m \\
        1 &   (-1)^m 
     \end{bmatrix}
     \begin{bmatrix} L_0 Y_{\ell}^{-m} \\ L_0 Y_{\ell}^{m} \end{bmatrix}
 = \frac{1}{\sqrt2}
     \begin{bmatrix} 
        i & - i ({-}1)^m \\
        1 &   ({-}1)^m 
     \end{bmatrix}
     \begin{bmatrix} -m\; Y_{\ell}^{-m} \\ m\; Y_{\ell}^{m} \end{bmatrix}
     \nonumber\\
  &= \frac{1}{\sqrt2}
     \begin{bmatrix} 
        -i & - i (-1)^m \\
        -1 &   (-1)^m 
     \end{bmatrix}
     \begin{bmatrix} m\; Y_{\ell}^{-m} \\ m\; Y_{\ell}^{m} \end{bmatrix}
   = \begin{bmatrix} 
       -i \; m ~\frac{1}{\sqrt2} \left( Y_{\ell}^{-m} + (-1)^m Y_{\ell}^{m} \right) \\
        i \; m ~\frac{i}{\sqrt2} \left( i ~Y_{\ell}^{-m} - i ~(-1)^m ~Y_{\ell}^{m} \right)
     \end{bmatrix} \nonumber\\
   &=  -i~
      \begin{bmatrix} 
        m ~Y_{\ell, m} \\
       -m ~Y_{\ell, -m}
     \end{bmatrix}
    = -i~
       \begin{bmatrix} 
          0 & m \\
         -m & 0 
     \end{bmatrix}
      \begin{bmatrix} 
        Y_{\ell, -m} \\
        Y_{\ell, m} 
     \end{bmatrix}
  \label{eqn:diff_lz_rsh_m}
\end{align}
The above derivation exploits the conversion formulas 
The composition of eq.~\eqref{eqn:diff_lz_rsh_0}-\eqref{eqn:diff_lz_rsh_m} results into the formula of $\bL_{\ell z}$, which is encoded by the recurrent definition
\begin{equation}
  \begin{gathered}
  \bL_{\ell z}^{[1]} =
    \begin{bmatrix}
     0 & 0 & 1 \\
     0 & 0 & 0 \\
    -1 & 0 & 0 
  \end{bmatrix}
  \quad
  \bL_{\ell z}^{[m]} =
  \left[
  \begin{array}{ccccc}
         &        &        &        &  m\\[2pt]
         & \multicolumn{3}{c}{\multirow{3}{*}{%
           \fbox{\rule[-6pt]{0pt}{18pt} $\bL_{\ell z}^{[m-1]}$ }}} &   \\[2pt]
         &        &        &        &        \\[2pt]
         &        &        &        &        \\[2pt]
   -m    &        &        &        &  
  \end{array}
  \right]
  \end{gathered}
  \label{eqn:lie_lz}
\end{equation}

The same procedure is applied in the case of operator $L_{+} = L_{x} + i L_{y}$.
The action of $L_{+}$ yields a linear combination 
We observe that $(L_{x} + i L_{y}))\bY_{\ell} = -i  (\bL_{\ell, x} + i \bL_{\ell, y}) \bY_{\ell} = (\bL_{\ell, y} - i \bL_{\ell, x}) \bY_{\ell}$. 
Thus, $\bL_{\ell x}$ and $\bL_{\ell y}$ are obtained respectively from the opposite of the imaginary part and the real part of the transformed RSH. 
First, the action on  $Y_{\ell, 0}$ yields 
\begin{align}
  L_{+} Y_{\ell, 0} 
  &= L_{+} Y_{\ell}^{0} 
   = \gamma_{\ell,0} Y_{\ell}^{1} 
   = -\frac{\gamma_{\ell,0} }{\sqrt{2}} \left( Y_{\ell, 1} + i ~Y_{\ell, -1} \right) \\
  &=  -\frac{\gamma_{\ell,0} }{\sqrt{2}} Y_{\ell, 1} - i  \frac{\gamma_{\ell,0} }{\sqrt{2}}  ~Y_{\ell, -1}
  \label{eqn:diff_lp_rsh_0}
\end{align}
As stated above the real part of  eq.~\eqref{eqn:diff_lp_rsh_0} defines the row $0$ of generator matrix $\bL_{\ell y}$ and the imaginary part the row $0$ of the generator matrix $\bL_{\ell x}$ after removing the multiplying factor $-i$. 
The same reasoning holds for the equations of operator $L_{+}$  when applied to RSH $Y_{\ell, -m}$ and $Y_{\ell, m}$.
First, we convert the RSH to CSH using the relation in appendix \ref{sec:appendix_csh_conversion}.
Second, we use the formulas of operator $L_{+}$ on CSH in eq.~\eqref{eqn:operator_action_to_csh}.
Finally, the result is converted to RSH again.
We start from the case $1 < m \leqslant \ell$:
\begin{align}
  L_{+} \begin{bmatrix} Y_{\ell, -m} \\ Y_{\ell, m} \end{bmatrix}
  &= \frac{1}{\sqrt2}
     \begin{bmatrix} i & - i (-1)^m \\ 1 & (-1)^m \end{bmatrix}
     \begin{bmatrix} L_{+} Y_{\ell}^{-m} \\ L_{+} Y_{\ell}^{m} \end{bmatrix} \\
&= \frac{1}{\sqrt2}
     \begin{bmatrix} i & - i (-1)^m \\ 1 & (-1)^m \end{bmatrix}
     \begin{bmatrix} \gamma_{\ell, -m} Y_{\ell}^{-m+1} \\ \gamma_{\ell, m} Y_{\ell}^{m+1} \end{bmatrix} \nonumber\\
  &= \frac{1}{\sqrt2}
     \begin{bmatrix} i & - i (-1)^m \\ 1 & (-1)^m \end{bmatrix}
     \begin{bmatrix}
       \frac{\gamma_{\ell, -m}}{\sqrt2} \left( Y_{\ell, m-1} - i\, Y_{\ell, -(m-1)} \right) \\[4pt]
       \frac{(-1)^{m+1}\gamma_{\ell, m}}{\sqrt2} \left( Y_{\ell, m+1} + i\, Y_{\ell, -(m+1)} \right)
     \end{bmatrix} \nonumber\\
  &= \frac{1}{2}
     \begin{bmatrix}
        \gamma_{\ell, -m} Y_{\ell, -(m-1)} - \gamma_{\ell, m} Y_{\ell, -(m+1)} \\
        \gamma_{\ell, -m} Y_{\ell, m-1}    - \gamma_{\ell, m} Y_{\ell, m+1}
     \end{bmatrix} \nonumber\\
  &\quad + \frac{i}{2}
     \begin{bmatrix}
        \gamma_{\ell, -m} Y_{\ell, m-1} + \gamma_{\ell, m} Y_{\ell, m+1} \\
       -\gamma_{\ell, -m} Y_{\ell, -(m-1)} - \gamma_{\ell, m} Y_{\ell, -(m+1)}
     \end{bmatrix}.
  \label{eqn:diff_lp_rsh_m}
\end{align}
Overall the action of operators $L_{x}$ and $L_{y}$ in the general case is
\begin{align}
  L_x Y_{\ell,-m} &= -\tfrac12\left(\gamma_{\ell,-m} Y_{\ell, m-1} + \gamma_{\ell,m} Y_{\ell, m+1}\right) \\
  L_x Y_{\ell, m} &= \tfrac12\left(\gamma_{\ell,-m} Y_{\ell,-(m-1)} + \gamma_{\ell,m} Y_{\ell,-(m+1)}\right) \\
  L_y Y_{\ell,-m} &= \tfrac12\left(\gamma_{\ell,-m} Y_{\ell,-(m-1)} - \gamma_{\ell,m} Y_{\ell,-(m+1)}\right) \\
  L_y Y_{\ell, m} &= \tfrac12\left(\gamma_{\ell,-m} Y_{\ell, m-1} - \gamma_{\ell,m} Y_{\ell, m+1}\right)
\end{align}
The case $m = 1$ must be handled apart, since $Y_{\ell, 0}$ has specific conversion equations:
\begin{align}
  L_{+} \begin{bmatrix} Y_{\ell, -1} \\ Y_{\ell, 1} \end{bmatrix}
  &= \frac{1}{\sqrt2}
     \begin{bmatrix} i & i \\ 1 & -1 \end{bmatrix}
     \begin{bmatrix} \gamma_{\ell,-1} Y_{\ell,0} \\ \gamma_{\ell,1}\, \frac{1}{\sqrt2}\left( Y_{\ell,2} + i Y_{\ell,-2} \right) \end{bmatrix} \nonumber\\
  &= \begin{bmatrix}
       -\frac{\gamma_{\ell,1}}{2} Y_{\ell,-2} \\
        \frac{\gamma_{\ell,-1}}{\sqrt2} Y_{\ell,0} - \frac{\gamma_{\ell,1}}{2} Y_{\ell,2}
     \end{bmatrix}
   + i \begin{bmatrix}
        \frac{\gamma_{\ell,-1}}{\sqrt2} Y_{\ell,0} + \frac{\gamma_{\ell,1}}{2} Y_{\ell,2} \\
       -\frac{\gamma_{\ell,1}}{2} Y_{\ell,-2}
     \end{bmatrix}.
  \label{eqn:diff_lp_rsh_1}
\end{align}
The coefficient is substituted as $\gamma_{\ell,-1}=\gamma_{\ell,0}$. 
Finally, the generators $\bL_{\ell x}$ and $\bL_{\ell y}$ are obtained by rearranging the coefficients obtained from ladder operator into rows and columns as
\begin{equation}
  \begin{gathered}
  \bL_{\ell x}^{[1]} =
    \begin{bmatrix}
     0                               & -\frac{\gamma_{\ell,0}}{\sqrt2} & 0 \\
     \frac{\gamma_{\ell,0}}{\sqrt2}  & 0                  & 0 \\
     0                               & 0                  & 0 
  \end{bmatrix}
  \quad
  \bL_{\ell x}^{[m]} =
  \left[
  \begin{array}{ccccc}
                       &        &        & -\frac{\gamma_{\ell,m}}{2}       &  0\\[2pt]
                       & \multicolumn{3}{c}{\multirow{3}{*}{%
                    \fbox{\rule[-6pt]{0pt}{18pt} $\bL_{\ell x}^{[m-1]}$ }}} &   -\frac{\gamma_{\ell,m}}{2} \\[2pt]
                       &        &        &        &        \\[2pt]
     \frac{\gamma_{\ell,m}}{2}  &        &        &        &        \\[2pt]
     0                 &  \frac{\gamma_{\ell,m}}{2}  &        &        &  
  \end{array}
  \right]
  \end{gathered}
  \label{eqn:lie_lx}
\end{equation}
\begin{equation}
  \begin{gathered}
  \bL_{\ell y}^{[1]} =
    \begin{bmatrix}
     0  & 0 & 0 \\
     0  & 0                  & -\frac{\gamma_{\ell,0}}{\sqrt2} \\
     0  & \frac{\gamma_{\ell,0}}{\sqrt2}                  & 0 
  \end{bmatrix}
  \quad
  \bL_{\ell y}^{[m]} =
  \left[
  \begin{array}{ccccc}
        0               &  \frac{\gamma_{\ell,m}}{2}      &        &       &  \\[2pt]
       -\frac{\gamma_{\ell,m}}{2}                & \multicolumn{3}{c}{\multirow{3}{*}{%
                    \fbox{\rule[-6pt]{0pt}{18pt} $\bL_{\ell y}^{[m-1]}$ }}} &    \\[2pt]
                       &        &        &        &        \\[2pt]
       &        &        &        &  -\frac{\gamma_{\ell,m}}{2} \\[2pt]
                      &    &        & \frac{\gamma_{\ell,m}}{2} &  0
  \end{array}
  \right]
  \end{gathered}
  \label{eqn:lie_ly}
\end{equation}
Given the three Cartesian generators $\bL_{\ell,x}$, $\bL_{\ell,y}$, and $\bL_{\ell,z}$, the Lie algebra representation $\bL_{\ell,\ba}$ of rotation $\ba$ is obtained as in eq.~\eqref{eqn:so3_band_ell_element}. 
Then, the RSH-rot matrix for band $\ell$ is computed through exponential mapping $\bM_{\ell, \ba} = \EXP{\bL_{\ell, \ba}}$. 
This is the first key result obtained through Lie algebra approach presented in this paper. 

There is an established literature about computation of matrix exponential and, for practical purpose where the band $\ell$ is limited, such approach is certainly robust and convenient. 
However, any method for exponential operation suffers from some kind of numerical weakness~\cite{MolerVanLoan03}.
In the remainder of this paper, we investigate the properties of Lie algebra representation of RSH-rot and formulate an alternative computation algorithm.

\subsection{Properties of Lie brackets}
\label{sec:rsh_lie_algebra_properties}

The Lie bracket is the constituting operation of a Lie algebra. 
By definition the representation $\bL_{\ell}$ preserves the brackets in the set of $(2 \ell+1)$-dimension skew symmetric matrices as in eq.~\eqref{eqn:lie_bracket_representation}. 
Let $\ba, \bb \in \mathbb{R}^3$ be two Rodrigues vectors, each associated with a rotation, and $\bL_{\ell,\ba}$, $\bL_{\ell,\bb}$ be their Lie algebra representations for band $\ell$. 
In this section, we present some properties related to the action of $\bL_{\ell,\ba}$ and of its group on $\bL_{\ell,\bb}$ (and vice versa). 
The \emph{adjoint representation} of the Lie algebra  is defined as
\begin{equation}
  \ADJ{\ba}{}(\bL_{\ell, \bb}) = [\bL_{\ell, \ba}, \bL_{\ell, \bb}]. 
\end{equation}
To avoid cumbersome notation, the subscript consists of the Rodrigues vector $\ba$ instead of the usual $\bL_{\ell, \ba}$. 
Due to the properties of Lie brackets, the adjoint is linear. 
Multiple compositions of such operator are denoted as 
\begin{equation}
  \ADJ{\ba}{n}(\bL_{\ell,\bb}) = [\underbrace{ \bL_{\ell,\ba}, [\dots, [\bL_{\ell,\ba}}_{n~times}, \bL_{\ell,\bb}] \ldots ]]
\end{equation}
For consistency, $\ADJ{\ba}{0}(\bL_{\ell,\bb}) = \bL_{\ell,\bb}$. 
In the following, we adopt the Cartesian basis of Lie algebra with the generators derived in section~\ref{sec:rsh_lie_algebra_generators}. 
Under such base there is a nice geometric tie between the Lie brackets and the cross product of Rodrigues vectors.

\begin{prop}
  Let $\bL_{\ell,\ba}, \bL_{\ell,\bb} \in \mathfrak{so}(2\ell+1)$ be two RSH-rot representations, then
  \begin{equation}
  	[\bL_{\ell,\ba}, \bL_{\ell,\bb}] = \bL_{\ell,\ba \times \bb}
  \end{equation}
  \label{thm:adjoint_cross}
\end{prop} 
\begin{proof}
  It is sufficient to expand the expression, to use the equations between the generators and to collect the terms:
  \begin{align}
     [\bL_{\ell,\ba}, \bL_{\ell,\bb}]
   & = [\sum_{i=1}^3 a_i ~\bL_{l,i}, \sum_{j=1}^3 b_j ~\bL_{l,j}] 
     = \sum_{i=1}^3 \sum_{j=1}^3 a_i~  b_j ~[ \bL_{l,i}, \bL_{l,j}] \\
   & = \sum_{i=1}^3 \sum_{j=i+1}^3 (a_i ~b_j - b_j ~ a_i) [ \bL_{l,i}, \bL_{l,j} ] \\
   & = (a_y b_z - b_z a_y) \bL_{l,x} + (a_z b_x - b_x a_z) \bL_{l,y} + (a_x b_y - b_y a_x)  \bL_{l,z} \\
   & = \bL_{\ell,\ba \times \bb}
  \end{align}
  This proves the relation between Lie bracket and cross product. 
\end{proof}

\begin{prop}
   Let $\bL_{\ell,\ba}, \bL_{\ell,\bb} \in \mathfrak{so}(2\ell+1)$ be two RSH-rot representations, then
  \begin{align}
 	\ADJ{\ba}{2}(\bL_{\ell, \bb}) &= \bL_{\ell, (\ba\transpose \bb) \ba - \| \ba \|^2 \bb} \\
 	\ADJ{\ba}{3}(\bL_{\ell, \bb}) &= -\| \ba \|^2 \ADJ{\ba}{}(\bL_{\ell, \bb}) 
  \end{align}
\end{prop} 
\begin{proof}
  The first part is the direct application of Prop.~\ref{thm:adjoint_cross} and cross product properties:
  \begin{align}
   \ADJ{\ba}{2}(\bL_{\ell, \bb})
    &= [\bL_{\ell, \ba}, [\bL_{\ell, \ba}, \bL_{\ell, \bb}]]
     = \bL_{\ell, \ba \times (\ba \times \bb) } \nonumber\\
    &= \bL_{\ell, (\ba\transpose \bb) \ba - \| \ba \|^2 \bb}
  \end{align}
  The second relation is obtained by recursively applying the same operation: 
  \begin{align}
  	\ADJ{\ba}{3}(\bL_{\ell, \bb})
  	&= [\bL_{\ell, \ba}, \ADJ{\ba}{2}(\bL_{\ell, \bb})]
  	= \bL_{\ell. (\ba \times ((\ba\transpose \bb) \ba - \| \ba \|^2 \bb))} \nonumber\\
  	&= \bL_{\ell, 0 - \| \ba \|^2 (\ba \times \bb)}
  	 = -\| \ba \|^2 \ADJ{\ba}{}(\bL_{\ell, \bb}) 
  \end{align}
  This proves the second formula. 
\end{proof}

The following corollary is even more useful.
It gives a direct geometric interpretation with the norm of Rodrigues vector $\theta_a$. 
Moreover, it traces back all the odd and even powers of the adjoint respectively to orders $1$ and $2$. 

\begin{corollary}
  Let $\ba = \theta_a \hat{\ba}$ and $\bb = \theta_b \hat{\bb}$ be two Rodrigues vectors, $\bL_{\ell,\ba}$ and $\bL_{\ell,\bb} \in \mathfrak{so}(2\ell+1)$ be their corresponding RSH-rot representations, then for every $k \geqslant 0$
  \begin{align}
  	\ADJ{\ba}{2 k+1}(\bL_{\ell, \bb}) &= (- \theta_a^2)^k ~\ADJ{\ba}{}(\bL_{\ell, \bb}) \\
  	\ADJ{\ba}{2 k+2}(\bL_{\ell, \bb}) &= (-\theta_a^2)^k ~\ADJ{\ba}{2}(\bL_{\ell, \bb}) 
  \end{align}
  \label{thm:adjoint_power_odd_even}
\end{corollary}

The following proposition provides an alternative expansion of the power of the adjoint. 

\begin{prop}
  Let $\bL_{\ell, \ba} = \ba \cdot \bL_{\ell}$ and $\bL_{\ell, \bb} = \bb \cdot \bL_{\ell}$, then
  \begin{equation}
	\ADJ{\ba}{n} (\bL_{\ell, \bb})
	= \sum_{k = 0}^n (-1)^k \binom{n}{k} \bL_{\ell, \ba}^{n-k} \bL_{\ell, \bb} \bL_{\ell, \ba}^k 
  \end{equation}
  \label{thm:power_adjoint}
\end{prop} 
\begin{proof}
  It can be easily verified for $n=1$ through the commutator definition $[\bL_{\ell, \ba},\bL_{\ell, \bb}] = \bL_{\ell, \ba} \bL_{\ell, \bb} - \bL_{\ell, \bb} \bL_{\ell, \ba}$. 
  We assume that the equation holds for $n$. 
  \begin{align}
  & \ADJ{\ba}{n+1} (\bL_{\ell, \bb})
    = [\bL_{\ell, \ba}, \ADJ{\ba}{n}(\bL_{\ell, \bb}) ] 
    = \left[\bL_{\ell, \ba}, \sum_{k = 0}^n (-1)^k \binom{n}{k} \bL_{\ell, \ba}^{n-k} \bL_{\ell, \bb} \bL_{\ell, \ba}^k \right] \\
  & = \sum_{k = 0}^n (-1)^k \binom{n}{k} \bL_{\ell, \ba}^{n-k+1} \bL_{\ell, \bb} \bL_{\ell, \ba}^k
     -\sum_{k = 0}^n (-1)^k \binom{n}{k} \bL_{\ell, \ba}^{n-k} \bL_{\ell, \bb} \bL_{\ell, \ba}^{k+1}\\
  & = \bL_{\ell, \ba}^{n+1} \bL_{\ell, \bb} + \sum_{k = 1}^n (-1)^k \binom{n}{k} \bL_{\ell, \ba}^{(n+1)-k} \bL_{\ell, \bb} \bL_{\ell, \ba}^k \nonumber \\
  &  + \sum_{k = 0}^{n-1} (-1)^{k+1} \binom{n}{k} \bL_{\ell, \ba}^{(n+1)-(k+1)} \bL_{\ell, \bb} \bL_{\ell, \ba}^{k+1}
     + (-1)^{n+1} \bL_{\ell, \bb} \bL_{\ell, \ba}^{n+1}
  \end{align}
  We substitute $k+1 \to k$ in the second summation.
  \begin{align}
  & = \bL_{\ell,\ba}^{n{+}1} \bL_{\ell, \bb}
    {+}\sum_{k{=}1}^n (-1)^k 
      \left(\binom{n}{k}{+}\binom{n}{k{-}1}\right) \bL_{\ell, \ba}^{(n{+}1)-k} \bL_{\ell, \bb} \bL_{\ell, \ba}^k
    {+}({-}1)^{n{+}1} \bL_{\ell, \bb} \bL_{\ell, \ba}^{n{+}1} \nonumber\\
  & = \bL_{\ell, \ba}^{n+1} \bL_{\ell, \bb}
  + \sum_{k = 1}^n (-1)^k \binom{n+1}{k} \bL_{\ell, \ba}^{(n+1)-k} \bL_{\ell, \bb} \bL_{\ell, \ba}^k
  + (-1)^{n+1} \bL_{\ell, \bb} \bL_{\ell, \ba}^{n+1} \nonumber\\
  &= \sum_{k=0}^{n+1} (-1)^k \binom{n+1}{k} \bL_{\ell, \ba}^{(n+1)-k} \bL_{\ell, \bb} \bL_{\ell, \ba}^k
  \end{align}
  This proves the proposition by induction. 
\end{proof}

Corollary~\ref{thm:adjoint_power_odd_even} and Proposition~\ref{thm:power_adjoint} can be combined together to obtain a Lie algebra version of the well known Rodrigues's rotation formula. 

\begin{prop}
  Let $\ba = \theta_a \hat{\ba}$ and $\bb = \theta_b \hat{\bb}$ be two Rodrigues vectors, $\bL_{\ell,\ba}$ and $\bL_{\ell,\bb} \in \mathfrak{so}(2\ell+1)$ be their corresponding RSH-rot representations, then
  \begin{align}
	\EXPSHORT{\bL_{\ell, \ba}} & \bL_{\ell, \bb} \EXPSHORT{-\bL_{\ell, \ba}} 
	= \sum_{n=0}^{\infty} \frac{1}{n!} \ADJ{\ba}{n}(\bL_{\ell, \bb}) \nonumber\\
    &= \bL_{\ell, \bb} + \frac{\sin \theta_a}{\theta_a} [\bL_{\ell, \ba}, \bL_{\ell, \bb}] 
	 + \frac{1 - \cos \theta_a}{\theta_a^2} [\bL_{\ell, \ba}, [\bL_{\ell, \ba}, \bL_{\ell, \bb}]]
     = \bL_{\ell, \bR_{\ba} \bb}
     \label{eqn:lie_conjugation_rodrigues}
  \end{align}
  \label{thm:conjugate_lie_alg}
\end{prop} 
\begin{proof}
  The proof starts from the expansion of exponential as 
  \begin{align}
  & \EXPSHORT{\bL_{\ell, \ba}} \bL_{\ell, \bb} \EXPSHORT{-\bL_{\ell, \ba}} 
    = \left( \sum_{i \geqslant 0} \frac{(\bL_{\ell, \ba})^i}{i!} \right) ~\bL_{\ell, \bb}~	
      \left( \sum_{j \geqslant 0} \frac{(-\bL_{\ell, \ba})^j}{j!} \right)	\\
  &= \sum_{i \geqslant 0} \sum_{j \geqslant 0} \frac{(-1)^j}{i! j!} \bL_{\ell, \ba}^i \bL_{\ell, \bb} \bL_{\ell, \ba}^j
  \end{align}
  We substitute $n=i+j$ and rename $k=j$ for readability:
  \begin{align}
  &= \sum_{n \geqslant 0} \frac{1}{n!} \sum_{k = 0}^n (-1)^k \binom{n}{k} \bL_{\ell, \ba}^{n-k} \bL_{\ell, \bb} \bL_{\ell, \ba}^k  
   = \sum_{n \geqslant 0} \frac{1}{n!} \ADJ{\ba}{n}(\bL_{\ell, \bb})
  \end{align}
  This step demonstrates the first equation. 
  Using Corollary~\ref{thm:adjoint_power_odd_even} we can collect the odd and even terms:
  \begin{align}
  & \EXPSHORT{\bL_{\ell, \ba}} \bL_{\ell, \bb} \EXPSHORT{-\bL_{\ell, \ba}} 
    = \bL_{\ell, \bb} + \sum_{k \geqslant 0} \frac{\ADJ{\ba}{2k+1}(\bL_{\ell, \bb})}{(2k+1)!} 
          + \sum_{k \geqslant 0} \frac{\ADJ{\ba}{2k+2}(\bL_{\ell, \bb})}{(2k+2)!}
  	\nonumber\\
  & = \bL_{\ell, \bb} + \sum_{k \geqslant 0} \frac{(-1)^k ~\theta_a^{2 k}}{(2k+1)!} \ADJ{\ba}{}(\bL_{\ell, \bb}) 
          + \sum_{k \geqslant 0} \frac{(-1)^k ~\theta_a^{2 k}}{(2k+2)!} \ADJ{\ba}{2}(\bL_{\ell, \bb}) 
    \nonumber\\
  & = \bL_{\ell, \bb} + \frac{ \ADJ{\ba}{}(\bL_{\ell, \bb}) }{\theta_a} 
       \underbrace{ \sum_{k \geqslant 0} \frac{(-1)^k ~\theta_a^{2 k+1}}{(2k+1)!} }_{\sin \theta_a}
          + \frac{ \ADJ{\ba}{2}(\bL_{\ell, \bb}) }{\theta_a^2} 
       \underbrace{ \sum_{k \geqslant 0} \frac{(-1)^k ~\theta_a^{2 k+2}}{(2k+2)!} }_{1 - \cos \theta_a}
  \end{align}
  The remaining direct substitution leads to the final formula. 
  Such Rodrigues-like formula is equivalent to the rotation of vector $\bb$ by the rotation corresponding to $\ba$, i.e. $\bR_{\ba}$. 
\end{proof}

Proposition~\ref{thm:conjugate_lie_alg} illustrates that the conjugation of a Lie algebra element $\bL_{\ell, \bb}$ by a Lie group element $\EXP{\bL_{\ell, \ba}}$ can be computed through operations on the Rodrigues vectors $\ba$ and $\bb$.
In the following, the formula is used to transform matrices into anti-diagonal form.  

\subsection{Computation of RSH Rotation through its Lie Algebra}

It is not trivial to derive a general formula to compute the exponential mapping of a Lie algebra representation matrix $\bL_{\ell,\ba}$. 
The only simple case is the rotation around axis $\hat{\bz}$. 
The derivation of RSH-rot matrix $\bM_{\ell,z}(\theta_z)$ has been discussed in section~\ref{sec:rsh_rot_z}. 
Now we show that $\bM_{\ell,z}(\theta_z)$ is the exponential of $\theta_z ~\bL_{\ell,z}$.
The generator $\bL_{\ell,z}$ is the anti-diagonal and skew-symmetric matrix in eq.~\eqref{eqn:lie_lz} and can be written as
\begin{equation}
  \bL_{\ell,z} 
  = \sum_{m=1}^{\ell} m ~\bB_{\ell,m}
  \qquad
  \left( \bB_{\ell,m} \right)_{rc} = \delta_{r,-m} \delta_{c,m} - \delta_{r,m} \delta_{c,-m}
\end{equation}
where $\left( \bB_{\ell,m} \right)_{rc}$ is the matrix element in row $r$ and column $c$ with $-\ell \leqslant r, c  \leqslant \ell$. 
The matrices $\bB_{\ell,m}$ are $(2\ell+1)$-dimensional, skew-symmetric and anti-diagonal with only two non-zero elements, $\left( \bB_{\ell,m} \right)_{-m,m} = 1$ and $\left( \bB_{\ell,m} \right)_{m,-m} = -1$.
The powers of $\bB_{\ell,m}^k$ have recurrent values
\begin{align}
  \begin{split}
  &\left( \bB_{\ell,m}^2 \right)_{rc} 
  = -\delta_{r,-m} \delta_{c,-m} - \delta_{r,m} \delta_{c,m}  \\
  &\bB_{\ell,m}^3 = -\bB_{\ell,m},
  \quad 
  \bB_{\ell,m}^4 = -\bB_{\ell,m}^2,
  \quad 
  \bB_{\ell,m}^5 = \bB_{\ell,m}, 
  \ldots
  \end{split}
\end{align}
They also satisfy the relation $\bB_{\ell,m} \bB_{\ell,m'} = \boldzero$ for $m' \neq m$. 
Thus, the $k$-th power of $\theta_z \bL_{\ell,z}$ for $k>0$ is equal to
\begin{align}
  (\theta_z \bL_{\ell,z})^k
  &= \sum_{m=1}^{\ell} (m \theta_z)^k ~\bB_{\ell,m}^k \\
  &= \sum_{m=1}^{\ell} (m \theta_z)^k ~(-1)^{\lfloor (k{-}1) / 2 \rfloor} ~\bB_{\ell,m}^{2 - (k\;\mathrm{mod}\;2)}
\end{align}
In the case $k=0$, the power is equal to the identity. 
The exponential map has the simple form related to the Rodrigues formula
\begin{align}
	& \EXP{\theta_z \bL_{\ell,z}}
	 = \bI + \sum_{k=1}^\infty \frac{ (\theta_z \bL_{\ell,z})^k }{ k! } \nonumber\\
	&= \bI + \sum_{m=1}^\ell \left(
	   \bB_{\ell,m} \sum_{k=0}^{\infty} (-1)^k \frac{(m \theta_z)^{2k+1}}{(2k+1)!}
	  {-}\bB_{\ell,m}^2 \sum_{k=1}^{\infty} (-1)^k \frac{(m \theta_z)^{2k}}{(2k)!} \right)\nonumber\\
	&= \bI + \sum_{m=1}^\ell \SIN{m \theta_z} \bB_{\ell,m} + \sum_{m=1}^\ell (1-\COS{m \theta_z}) ~\bB_{\ell,m}^2
\end{align}
After expanding all the above matrix terms, the outcome is equal to eq.~\eqref{eqn:rsh_rot_z}. 
We have verified the expected result. 

One possible approach to compute the exponential map of the representation $\bL_{\ell,\ba}$ of the algebra for the rotation $\bR_{\ba}$ is to transform it into an anti-diagonal matrix via a similarity relation. 
The similarity is given by an invertible matrix $\bF_{\ell, \hat{\ba}}$ such that
\begin{equation}
  \bF_{\ell,\hat{\ba}}^{-1} ~\bL_{\ell, \ba} ~\bF_{\ell,\hat{\ba}} = \theta_a ~\bL_{\ell,z} 
  \qquad
  \bL_{\ell, \ba} = \bF_{\ell,\hat{\ba}} \;  (\theta_a ~\bL_{\ell,z} ) \; \bF_{\ell,\hat{\ba}}\transpose
  \label{eqn:lie_algebra_antidiagonal}
\end{equation}
We call the second equation of eq.~\eqref{eqn:lie_algebra_antidiagonal} \emph{anti-diagonal factorization} of $\bL_{\ell, \ba}$, since the anti-diagonal matrix $\bL_{\ell,z}$ appears on the right side. 
Proposition~\ref{thm:conjugate_lie_alg} provides an easy way to find such matrix $\bF_{\ell,\hat{\ba}}$. 
Indeed, eq.~\eqref{eqn:lie_algebra_antidiagonal} is a conjugation by $\bF_{\ell,\hat{\ba}}^{-1} = \EXPSHORT{\bL_{\ell,\bc}}$, where $\bc$ is the Rodrigues vector $\bc$ of the rotation $\bR_{\bc}$ that transforms $\ba$ into $\theta_a \hat{\bz}$.
The vector $\bc = \theta_c \hat{\bc}$ is easily obtained from elementary geometric reasoning as
\begin{equation}
  \hat{\bc} = \frac{\ba \times \hat{\bz}}{\lVert \ba \times \hat{\bz} \rVert}
  \qquad
  \theta_c = \mathrm{acos}(\hat{\ba}\transpose \hat{\bz})
\end{equation}
(for $\hat{\ba}=\pm\hat{\bz}$ any axis orthogonal to $\hat{\bz}$ can be used for $\hat{\bc}$). 
It is clear that $\bF_{\ell,\hat{\ba}} = \EXPSHORT{-\bL_{\ell,\bc}}$ is orthogonal. 
Finally, eq.~\eqref{eqn:lie_algebra_antidiagonal} can be rearranged to write the exponential 
\begin{equation}
   \EXP{\bL_{\ell,\ba}}
   = \sum_{k=0}^{\infty} \frac{( \bF_{\ell,\hat{\ba}} (\theta_a \bL_{\ell,z})  \bF_{\ell,\hat{\ba}}\transpose)^k}{k!}
   = \bF_{\ell,\hat{\ba}} \; \bM_{\ell,z}(\theta_a) \; \bF_{\ell,\hat{\ba}}\transpose
\end{equation}
In this new equation the exponential is obtained from a product consisting of $\bF_{\ell, \hat{\ba}}$ and its inverse, and the RSH-rot $\bM_{\ell,z}(\theta_a)$. 
All these terms can be computed. 
Unfortunately, the advantage of such factorization is limited.
Since the goal is to avoid a general-purpose matrix exponential, and $\bF_{\ell, \hat{\ba}}$, although now fully characterized, is itself defined through the exponential map, the factorization is not yet a computational alternative. 
The next section is devoted to derivation of an algorithm to compute this factor $\bF_{\ell, \hat{\ba}}$.

\section{Anti-diagonal Factorization of RSH Rotation Matrix}

The previous section introduced the Lie algebra representation $\bL_{\ell}$ of the rotation action on SH and showed how it can be used to compute the RSH rotation matrix $\bM_{\ell}$ through the exponential map. Although this approach provides a satisfactory mean of evaluating the rotation matrix, an alternative formulation can be obtained by exploiting the common eigenvectors of $\bM_{\ell}$ and $\bL_{\ell}$. These eigenvectors allow both matrices to be factorized into simpler components. In this section, we derive a recurrence relation for computing the eigenvectors and use them to obtain a convenient factorization of the RSH rotation matrix.

\subsection{Eigenvectors of RSH Rotation Matrix and Anti-diagonal Factorization}
\label{sec:rsh_rot_eigenvectors}

The spectral structures of an RSH rotation matrix and its corresponding Lie algebra representation are closely related. This relationship makes it possible to perform the analysis in the Lie algebra domain, whose structure is comparatively simpler.

We adopt the notation previously introduced and represent a rotation equivalently by the matrix $\bR_\ba \in SO(3)$ and its Rodrigues vector $\ba = \theta_a \hat{\ba}$, where $\|\hat{\ba}\|=1$. The corresponding RSH rotation matrix for band $\ell$ is denoted by $\bM_{\ell}(\bR_\ba)$, while its Lie algebra representation is $\bL_{\ell,\ba}$, as defined in eq.~\eqref{eqn:so3_band_ell_element}. The exponential map relating these two representations leads to the following property.

\begin{prop}
  \label{thm:rsh_eigen}
  Let $\bM_{\ell}(\bR_\ba)$ and $\bL_{\ell,\ba}$ denote the representations of the rotation $\bR_\ba$ in the RSH domain, respectively for the rotation group and its Lie algebra. They share the same eigenvectors $\bv_m=\bt_m+i\bw_m\in\mathbb{C}^{2\ell+1}$, for $m=-\ell,\ldots,\ell$, satisfying
  \begin{equation}
    \bL_{\ell \ba} ~\bv_m = i ~m ~\theta_a ~\bv_m,
    \qquad
    \bM_\ell(\bR_\ba) ~\bv_m = \EXPSHORT{i ~m ~\theta_a} ~\bv_m.
    \label{eqn:eigenvector_rsh_lie_rot}
  \end{equation}
  Each eigenvalue $i\,m\,\theta_a$ of $\bL_{\ell,\ba}$ therefore corresponds to an eigenvalue $\EXPSHORT{i ~m ~\theta_a}$ of $\bM_{\ell}(\bR_\ba)$. Moreover, $\bv_{-m}=\bv_m^*$ for $0<m\leqslant\ell$.
\end{prop}

\begin{proof}
We define the diagonal matrix $\bLambda_\ell$ containing the eigenvalues and the matrix $\bV_\ell$ whose columns are the corresponding eigenvectors:
\begin{equation}
  \bLambda_\ell = \mathrm{diag}(-i \ell \theta_a, \dots, 0, \ldots, i \ell \theta_a)
  \qquad
  \bV_{\ell} = \bmat{\bv_{-\ell} & \ldots & \bv_{0} & \ldots & \bv_{\ell}}.
\end{equation}
By the definition of eigenvalues and eigenvectors, the Lie algebra element admits the factorization
$\bL_{\ell,\ba} = \bV_\ell \bLambda_\ell \bV_\ell\hermitian$.

The rotation matrix $\bM_{\ell}(\bR_\ba)$ is obtained by exponentiating $\bL_{\ell,\ba}$. Thus,
\begin{equation}
  \bM_\ell(\bR_\ba)
  = \EXP{\bL_{\ell,\ba}}
  = \sum_{k \geqslant 0} \frac{(\bV_\ell \bLambda_\ell \bV_\ell\hermitian)^k}{k!}
  = \bV_\ell\;\EXP{\bLambda_\ell}\;\bV_\ell\hermitian.
\end{equation}
The matrix $\EXP{\bLambda_\ell}$ is diagonal, with diagonal elements $\EXPSHORT{i\,m\,\theta_a}$ for $m=-\ell,\ldots,\ell$. These are precisely the eigenvalues of $\bM_\ell(\bR_\ba)$ derived in Proposition~\ref{thm:rsh_eigenvalues}.

Finally, since $\bL_{\ell,\ba}$ is a real matrix, its eigenvalues associated with $m$ and $-m$ are complex conjugates. The corresponding eigenvectors can therefore be chosen as conjugate pairs, such that
$\bv_{-m}=\bv_m^*=\bt_m-i\bw_m$.
\end{proof}

Proposition~\ref{thm:rsh_eigen} establishes that the group representation and its Lie algebra representation share the same eigenvectors. We can therefore focus on the normalized Lie algebra element
$\bL_{\ell,\hat{\ba}}=\bL_{\ell,\ba}/\theta_a$, which provides a simpler starting point for the analysis. Its eigendecomposition makes it possible to express the matrix in a diagonal basis.

\begin{prop}
  \label{thm:rsh_antidiagonal}
  Let $\bL_{\ell,\ba}$ be a Lie algebra representation associated with the Rodrigues vector $\ba=\theta_a\hat{\ba}$ and having eigenvectors $\bv_m=\bt_m+i\bw_m$. Then, it admits the factorization
  \begin{equation}
     \bL_{\ell,\ba}
     = \bF_{\ell}\;(\theta_a\bL_{\ell,z})\;\bF_{\ell}\transpose
     \label{eqn:anti_diagonal_factorization}
  \end{equation}
  where $\bL_{\ell,z}$ is the anti-diagonal generator associated with a rotation about the axis $\hat{\bz}$, as defined in eq.~\eqref{eqn:lie_lz}, and $\bF_{\ell}$ is given by
  \begin{equation}
    \bF_{\ell}
     = \begin{bmatrix}
         \sqrt{2}\,\bt_{\ell} & \cdots & \sqrt{2}\,\bt_{1} &
         \bt_{0} &
         \sqrt{2}\,\bw_{1} & \cdots & \sqrt{2}\,\bw_{\ell}
       \end{bmatrix}.
  \end{equation}
\end{prop}

\begin{proof}
The matrix can be decomposed in terms of its eigenvalues and eigenvectors, as established in Proposition~\ref{thm:rsh_eigen}:
\begin{equation}
  \bL_{\ell,\ba}
   = \sum_{m=-\ell}^{\ell} \lambda_m \;\bv_m\bv_m\hermitian
   = \sum_{m=1}^{\ell} i\,m\,\theta_a
     \left(
       \bv_m\bv_m\hermitian
       - \bv_{-m}\bv_{-m}\hermitian
     \right).
\end{equation}
The contribution associated with the zero eigenvalue vanishes and has therefore been omitted. The two conjugate eigenvectors can be combined as follows:
\begin{align}
   \bv_m\bv_m\hermitian
   &= (\bt_m\bt_m\transpose+\bw_m\bw_m\transpose)
    + i(\bw_m\bt_m\transpose-\bt_m\bw_m\transpose), \\
   \bv_{-m}\bv_{-m}\hermitian
   &= (\bt_m\bt_m\transpose+\bw_m\bw_m\transpose)
    - i(\bw_m\bt_m\transpose-\bt_m\bw_m\transpose).
\end{align}
Substituting these expressions into the eigendecomposition yields
\begin{align}
  \bL_{\ell,\ba}
   &= \sum_{m=1}^{\ell}
      2m\,\theta_a
      (\bt_m\bw_m\transpose-\bw_m\bt_m\transpose) \\
   &= \sum_{m=1}^{\ell}
      (m\theta_a)(\sqrt{2}\,\bt_m)(\sqrt{2}\,\bw_m)\transpose
      + \sum_{m=1}^{\ell}
      (-m\theta_a)(\sqrt{2}\,\bw_m)(\sqrt{2}\,\bt_m)\transpose \\
   &= \bF_{\ell}\;(\theta_a\bL_{\ell,z})\;\bF_{\ell}\transpose.
\end{align}
The eigenvectors are normalized such that $\|\bv_m\|=1$. Consequently, their real and imaginary parts satisfy
$\|\bt_m\|=\|\bw_m\|=1/\sqrt{2}$.
The factor of $2$ in the preceding expression is therefore distributed between the two vectors to obtain the normalization used in the definition of $\bF_\ell$.
\end{proof}

The anti-diagonal factorization is thus obtained from the spectral decomposition of $\bL_{\ell,\ba}$. Once $\bF_{\ell}$ has been computed, the corresponding RSH rotation matrix can be evaluated as
$\bM_{\ell}(\bR_\ba)=\bF_{\ell}\;\bM_{\ell,z}(\theta_a)\;\bF_{\ell}\transpose$.
In the following subsection, we introduce an alternative procedure for computing the eigenvectors based on a recurrence relation.

\subsection{Eigenvector Evaluation by Recurrence Relations}

To diagonalize $\bL_{\ell,\hat{\ba}}$, we introduce three mutually orthogonal axes
$\hat{\be}_1,\hat{\be}_2,\hat{\be}_3\in\mathbb{R}^3$, with
$\hat{\be}_3=\hat{\ba}$ aligned with the rotation axis and
$\hat{\be}_1,\hat{\be}_2$ orthogonal to it.

If $\hat{\ba}=\hat{\bz}$, the rotation axis coincides with the $z$-axis and
$\rho=\sqrt{\hat{a}_x^2+\hat{a}_y^2}=0$. In this case, the matrix
$\bM_{\ell}$ already has the simple form given in eq.~\eqref{eqn:rsh_rot_z}.

We therefore consider the more general case in which $\rho>0$. One possible choice for the axes $\hat{\be}_1$ and $\hat{\be}_2$ is
\begin{equation}
  \hat{\be}_1 =
  \frac{1}{\rho}
  \begin{bmatrix}
    -\hat{a}_y \\ \hat{a}_x \\ 0
  \end{bmatrix},
  \qquad
  \hat{\be}_2 =
  \frac{1}{\rho}
  \begin{bmatrix}
    -\hat{a}_x\hat{a}_z \\
    -\hat{a}_y\hat{a}_z \\
    \rho^2
  \end{bmatrix}.
  \label{eqn:orto_basis_e1_e2}
\end{equation}
It is straightforward to verify that these axes are mutually orthogonal and satisfy the right-handedness condition
$\hat{\be}_j\times\hat{\be}_k=\epsilon_{jkh}\hat{\be}_h$.
The corresponding Lie algebra representations are defined as
$\hat{\bL}_j=\bL_{\ell,\hat{\be}_j}$, for $j=1,2,3$.
Using Proposition~\ref{thm:adjoint_cross}, their Lie brackets satisfy
\begin{equation}
  [\hat{\bL}_j,\hat{\bL}_k]
  =(\hat{\be}_j\times\hat{\be}_k)\cdot\bL_\ell
  =(\epsilon_{jkh}\hat{\be}_h)\cdot\bL_\ell
  =\epsilon_{jkh}\hat{\bL}_h.
\end{equation}
The generators $\hat{\bL}_1$, $\hat{\bL}_2$, and $\hat{\bL}_3$ therefore form an alternative Cartesian basis. We now introduce the complex generators
$\hat{\bL}_{\pm}=\hat{\bL}_2\pm i\hat{\bL}_1$, replacing the first two generators with their complex combinations. 
$\hat{\bL}_{+}$ and $\hat{\bL}_{-}$ must not be confused with the ladder operators $L_\pm$ of eq.~\eqref{eqn:operator_definition} or with the generators $\bL_{\ell\pm}$ of the complex basis. 
The resulting basis satisfies the following Lie bracket relations:
\begin{equation}
  [\hat{\bL}_{3},\hat{\bL}_{+}]
  = i\hat{\bL}_{+},
  \quad
  [\hat{\bL}_{3},\hat{\bL}_{-}]
  = -i\hat{\bL}_{-},
  \quad
  [\hat{\bL}_{+},\hat{\bL}_{-}]
  = 2i\hat{\bL}_{3}.
  \label{eqn:rot_axis_lie_brackets}
\end{equation}
These relations show that the commutator of $\hat{\bL}_3$ with either $\hat{\bL}_+$ or $\hat{\bL}_-$ is proportional to that same generator. This property is the key to deriving a recurrence relation for the eigenvectors.

\begin{prop}
  \label{thm:rsh_eigenvector_recursive}
  Let $\bv_m$ be an eigenvector of $\hat{\bL}_3=\bL_{\ell,\hat{\ba}}$ with eigenvalue $i\,m$. Then, the eigenvectors $\bv_{m+1}$ and $\bv_{m-1}$, associated with the eigenvalues $i(m+1)$ and $i(m-1)$, respectively, satisfy
  \begin{equation}
    \bv_{m+1}\propto\hat{\bL}_{+}\bv_m,
    \qquad
    \bv_{m-1}\propto\hat{\bL}_{-}\bv_m.
    \label{eqn:rsh_eig_recur_propto}
  \end{equation}
\end{prop}

\begin{proof}
Assume that $\bv_m$ is an eigenvector of $\hat{\bL}_3$ with eigenvalue $i\,m$, i.e.,
$\hat{\bL}_3\bv_m=i\,m\,\bv_m$.
Applying the commutator relations in eq.~\eqref{eqn:rot_axis_lie_brackets} to $\bv_m$ gives
\begin{align}
  \pm i\hat{\bL}_{\pm}\bv_m
  &= [\hat{\bL}_3,\hat{\bL}_{\pm}]\bv_m \\
  &= \hat{\bL}_3\hat{\bL}_{\pm}\bv_m
     -\hat{\bL}_{\pm}\hat{\bL}_3\bv_m \\
  &= \hat{\bL}_3\hat{\bL}_{\pm}\bv_m
     -\hat{\bL}_{\pm}(i\,m\,\bv_m).
\end{align}
In the last equality, we have used the eigenvector equation for $\bv_m$. Rearranging the terms and collecting the contribution involving $\hat{\bL}_{\pm}\bv_m$ yields
\begin{align}
  \hat{\bL}_3(\hat{\bL}_{\pm}\bv_m)
  = i(m\pm1)(\hat{\bL}_{\pm}\bv_m).
\end{align}
Therefore, if $\bv_m$ is an eigenvector of $\hat{\bL}_3$ with eigenvalue $i\,m$, then $\hat{\bL}_{+}\bv_m$ and $\hat{\bL}_{-}\bv_m$ are, respectively, unnormalized eigenvectors associated with the eigenvalues $i(m+1)$ and $i(m-1)$.
\end{proof}

Proposition~\ref{thm:rsh_eigenvector_recursive} provides a procedure for recursively computing the eigenvectors. Starting from an initial eigenvector $\bv_0$, the subsequent eigenvectors can be obtained by repeatedly applying $\hat{\bL}_{+}$ and normalizing the result. The initial eigenvector is
$\bv_0=\bY_{\ell}(\hat{\ba})$, corresponding to the unique real eigenvalue $1$ of $\bM_{\ell}(\bR_\ba)$, as established in Proposition~\ref{thm:rsh_band_l_rot_eigenvector}, and to the eigenvalue $0$ of $\hat{\bL}_{3}$ from Proposition~\ref{thm:rsh_eigen}.
The procedure terminates once $\bv_{\ell}$ has been computed. 
Observe that the eigenvectors are not unique, but each $\bv_m$ is defined only up to an individual phase $\EXPSHORT{i\beta_m}$.
Indeed, there is one degree of freedom in the choice of $\hat{\be}_1$ and $\hat{\be}_2$ in eq.~\eqref{eqn:orto_basis_e1_e2}, which rotates each pair $(\bt_m, \bw_m)$ within its invariant plane.

The same recursive procedure can also be used to construct the factorization of the RSH rotation matrix. In Section~\ref{sec:rsh_rot_z}, we showed that a rotation matrix can be decomposed as
$\bR=\bF\;\bR_z\;\bF\transpose$, where $\bR_z$ represents a rotation about the $z$-axis. By exploiting the real and imaginary parts of the eigenvectors, we can construct a matrix $\bF_{\ell}(\hat{\ba})$ such that
\begin{equation}
  \bM_{\ell}(\bR_\ba)
  = \bF_{\ell}(\hat{\ba})\;
    \bM_{\ell z}(\theta_a)\;
    \bF_{\ell}(\hat{\ba})\transpose.
\end{equation}
The exponential-map relationship also gives the corresponding factorization of the Lie algebra representation:
\begin{equation}
  \bL_{\ell,\ba}
  = \bF_{\ell}(\hat{\ba})\;
    (\theta_a\bL_{\ell z})\;
    \bF_{\ell}(\hat{\ba})\transpose.
\end{equation}

To construct $\bF_{\ell}(\hat{\ba})$, we consider the real and imaginary parts of the eigenvectors $\bv_m=\bt_m+i\bw_m$. The eigenvector equation associated with the eigenvalue $i\,m$ yields the following relationships between these two components:
\begin{equation}
  \bt_m = \frac{1}{m}\bL_{\ell,\hat{\ba}}\bw_m,
  \qquad
  \bw_m = -\frac{1}{m}\bL_{\ell,\hat{\ba}}\bt_m.
\end{equation}
These equations constrain $\bt_m$ and $\bw_m$ and will be used in the following derivation.

The recurrence relation involving $\hat{\bL}_{+}\bv_m$, given in eq.~\eqref{eqn:rsh_eig_recur_propto}, can be normalized using the coefficient $\gamma_{\ell m}$ introduced in eq.~\eqref{eqn:operator_action_to_csh}. It can therefore be written as
\begin{align}
  \bv_{m+1}
  &= \bt_{m+1}+i\bw_{m+1} \\
  &= \frac{1}{\gamma_{\ell m}}
     (\hat{\bL}_{2}+i\hat{\bL}_1)
     (\bt_m+i\bw_m).
\end{align}
Separating the real and imaginary parts gives
\begin{align}
  \bt_{m+1}
  &= \frac{1}{\gamma_{\ell m}}
     (\hat{\bL}_{2}\bt_m-\hat{\bL}_1\bw_m), \\
  \bw_{m+1}
  &= \frac{1}{\gamma_{\ell m}}
     (\hat{\bL}_{2}\bw_m+\hat{\bL}_1\bt_m).
\end{align}

This recurrence relation provides a means of computing the eigenvectors $\bv_m$ for $m=0,\ldots,\ell$, starting from $\bt_0=\bY_\ell(\hat{\ba})/\lVert\bY_\ell(\hat{\ba})\rVert$ and $\bw_0=\boldzero$.
In exact arithmetic the computed eigenvectors are mutually orthogonal.
In floating point, however, the recurrence is unstable: the round-off components of $\bv_m$ along the other eigenvectors $\bv_k$ are multiplied at each step by the ratios $\gamma_{\ell,k}/\gamma_{\ell,m}$ of the ladder coefficients, which exceed one for $-m\leqslant k<m$.
Without countermeasures, the orthogonality residual $\lVert\bF_\ell^\top\bF_\ell-\bI\rVert_F$ grows from about $10^{-13}$ at $\ell=16$ to about $10^{-8}$ at $\ell=64$.
The amplified components lie in the span of the vectors already computed, so they are removed by re-orthogonalizing each new pair $(\bt_{m+1},\bw_{m+1})$ against all previous $\bt_k,\bw_k$, $k\leqslant m$ (classical Gram--Schmidt with two passes), and normalizing $\lVert\bv_{m+1}\rVert=1$.
This keeps $\bF_\ell$ orthogonal to round-off level at a cost of $\mathcal{O}(\ell^3)$ operations per band, the same order as the product $\bF_\ell\bM_{\ell z}(\theta_a)\bF_\ell^\top$.
Once these eigenvectors have been obtained, their real and imaginary parts are used to construct $\bF_{\ell}(\hat{\ba})$ as in eq.~\eqref{eqn:anti_diagonal_factorization} and, consequently, to evaluate the anti-diagonal factorization of the RSH rotation matrix.

\section{Differentiation of RSH Rotation}
\label{sec:differentiation}

In the previous section, we have proposed two methods for the computation of the RSH-rot matrix, which is the representation of the action of $SO(3)$ on RSH coefficients.
The first uses the exponential map of the Lie algebra representation matrix whereas the second exploits the anti-diagonal factorization. 
Both are consequences of treating the rotation group at the level of its Lie algebra, which is where the tangent space of $SO(3)$ lives and where differentiation becomes algebraic rather than numerical.
This yields an advantage that the formulations currently in the literature do not offer. 
Because a function of the rotation is built from an RSH series and a rotated RSH series, and because the Lie algebra characterizes the tangent space of the group, the derivative of that function with respect to the rotation is available in closed form. 
No finite-difference approximation of the rotation, and no re-expansion of the RSH series at each perturbation, is required.

This advantage is apparent by considering a relevant application example. 
Several problems in engineering involve the directional distribution of quantities such as antenna powers or light reflectance on surfaces.
Such distributions are modeled as function on the unit sphere $f:\mathbb{S}^2 \to \mathbb{R}$. 
A notable example of such functions is the \emph{von Mises--Fisher} (vMF) distribution. 
Given a unitary vector $\bmu \in \mathbb{S}^2$ and a concentration parameter $\kappa$, the vMF is 
\begin{equation}
  f(\bphat) = \frac{\kappa}{4 \pi \sinh(\kappa)} \EXP{\kappa \bmu\transpose \bphat}
  \label{eqn:vmf_distribution}
\end{equation}
The vMF can be expanded into RSH  as in eq.~\eqref{eqn:func_rsh_series} with coefficients
\begin{equation}
  \bc_{\ell} = \lambda(\kappa) \; \bY_{\ell}(\bmu) = \sqrt{ \frac{\pi \kappa}{2} } \; \frac{ \BesselI{\ell+1/2}{\kappa} }{ \sinh{\kappa} }\; \bY_{\ell}(\bmu)
  \label{eqn:vmf_rsh_coeffs}
\end{equation}
where $\BesselI{\ell+1/2}{\cdot}$ is the modified Bessel function of the first kind of fractional order $\ell{+}1{/}2$. 
Observe that the vMF is a symmetric distribution and falls in the case addressed in Proposition~\ref{thm:rsh_axis_symmetric}. 
Indeed, the point on sphere $\bmu$ is the symmetry axis and the coefficients in eq.~\eqref{eqn:vmf_rsh_coeffs} are in the form of the product of constant $\lambda(\kappa)$ and $\bY_{\ell}(\bmu)$. 
A rotation $\bR_{\alpha \bmu}$ about $\bmu$ by any angle $\alpha$ leaves the distribution unchanged (it is sufficient to substitute $\bR_{\alpha \bmu} \bmu = \bmu$). 
This ambiguity is removed by considering mixtures of vMF distributions, whose RSH series expansion is
\begin{align}
  f(\bphat) 
  &= \sum_{j=1}^{n} w_j \; \frac{\kappa_j}{4 \pi \sinh(\kappa_j)} \EXP{\kappa_j \bmu_j\transpose \bphat} \\
  &= \sum_{\ell \geqslant 0} \left(\sum_{j=1}^{n} w_j \lambda(\kappa_j) \bY_{\ell}(\bmu_j) \right)\transpose \;\bY_{\ell}(\bphat) \label{eq:vmf_mixture_rsh}
\end{align}
where $w_j$, $\bmu_j$ and $\kappa_j$ are respectively the weight, the direction, and the concentration parameters of the $j$-th mode of the mixture. 
With distinct $\bmu_j$ there is no symmetry. 
This distribution is a valuable benchmark and is therefore also used in the experiments. 

Several distributions can be measured from different viewpoints related by a rotation in space. 
Let $f_s, f_d: \mathbb{S}^2 \to \mathbb{R}$ be two functions on unit spheres representing the directional distributions of the same quantity observed in two different locations called respectively source and destination. 
Being defined on $\mathbb{S}^2$ both the functions can be expanded in RSH as
\begin{equation}
  f_s(\bphat) = \sum_{\ell \geqslant 0} \bs_\ell\transpose ~\bY_{\ell}(\bphat) 
  \qquad
  f_d(\bphat) = \sum_{\ell \geqslant 0} \bd_\ell\transpose ~\bY_{\ell}(\bphat) 
\end{equation}
The vectors $\bs_\ell, \bd_\ell \in \mathbb{R}^{2\ell+1}$ are the RSH coefficients in each band $\ell$.
We assume that the source and destination distributions represent the same quantities observed from source and destination poses and the source pose can be aligned to destination pose by an unknown rotation $\bR_\ba \in SO(3)$. 
Indeed, except for noise and occlusion, $f_d$ is a rotated copy of $f_s$. 
One way to estimate such rotation is through a correlation function measuring the similarity between a rotated $f_s$ and $f_d$. 
Such correlation can be defined as the integral product 
\begin{align}
  C(\bR_\ba) 
  &= \langle f_d(\cdot), f_s(\bR_\ba \cdot) \rangle 
   = \int_{\mathbb{S}^2} f_d(\bphat) \; f_s(\bR_\ba \bphat) \; d\bphat \\
  &= \int_{\mathbb{S}^2} 
     ~\left( \sum_{\ell \geqslant 0} \bd_\ell\transpose ~\bY_{\ell}(\bphat) \right)  
     ~\left( \sum_{\ell \geqslant 0} \bs_\ell\transpose ~\bY_{\ell}(\bR_\ba \bphat) \right) 
     ~d\bphat  \\
  &= \sum_{\ell \geqslant 0} \bd_\ell\transpose  
     \left( \int_{\mathbb{S}^2} \bY_{\ell}(\bR_\ba \bphat) \bY_{\ell}\transpose(\bphat) d\bphat \right) 
     \bs_\ell 
   = \sum_{\ell \geqslant 0} \bd_\ell\transpose ~\bM_\ell(\bR_\ba) ~\bs_\ell 
  \label{eqn:rsh_correlation}
\end{align}
In the last step we have used the orthogonality of RSH of different bands and substituted the definition~\eqref{eqn:rsh_rot_def} of $\bM_\ell$. 
The rotation corresponding to the best alignment is the maximum of the objective function $C(\bR_\ba)$. 

The proposed Lie algebra approach supports the estimation of the optimum in several ways. 
First, it allows the vector parameter $\ba$ to be used as the domain of $C(\cdot)$ instead of the rotation $\bR_\ba$ on the manifold $SO(3)$. 
Second, it provides the local expansion of the objective function required by most solvers. 
Given a point $\tilde{\ba}$ or $\bR_{ \tilde{\ba} }$, the increments on $\bM_\ell$ with respect to $\tilde{\ba}$ can be expressed by $\Delta \ba$ as
\begin{equation}
  \bM_\ell( \bR_{\tilde{\ba}} \; \EXPSHORT{ \bOmega_{\Delta\ba} } ) 
  = \bM_\ell( \bR_{\tilde{\ba}} ) ~\EXPSHORT{ \bL_{\ell, \Delta\ba} }
  = \bM_\ell( \bR_{\tilde{\ba}} ) 
    \left( \bI + \bL_{\ell, \Delta\ba} + \frac{ \bL_{\ell, \Delta\ba}^2 }{2} + \ldots\right)
  \label{eqn:rsh_exp_increment}
\end{equation}
Second-order iterative methods, such as Newton's method, update the estimate using a quadratic approximation of the objective function. 
The approximation of $C(\bR_\ba)$ is obtained from the substitution of eq.~\eqref{eqn:rsh_exp_increment} as
\begin{align}
  C( \bR_{\tilde{\ba}} \; \EXPSHORT{ \bOmega_{\Delta\ba} } ) 
  &= \sum_{\ell \geqslant 0} \bd_\ell\transpose 
    ~\bM_\ell( \bR_{\tilde{\ba}} ) 
    \left( \bI + \bL_{\ell, \Delta\ba} + \frac{ \bL_{\ell, \Delta\ba}^2 }{2} + \ldots\right)
    ~\bs_\ell \\
  &= c + \bb\transpose ~\Delta\ba + \frac{1}{2} \Delta\ba\transpose ~\bH ~\Delta\ba + \mathcal{O}(\| \Delta\ba \|^3)
  \label{eqn:corr_quadraitc_approx}
\end{align}
The three coefficients in the above quadratic function are equal to 
\allowdisplaybreaks
\begin{align}
  c 
  &= \sum_{\ell \geqslant 0} \tilde{\bd}_\ell\transpose ~\bs_\ell \\
  \bb 
  &= \begin{bmatrix}
  	   \sum_{\ell \geqslant 0} \tilde{\bd}_\ell\transpose ~\bL_{\ell x} ~\bs_\ell \\
  	   \sum_{\ell \geqslant 0} \tilde{\bd}_\ell\transpose ~\bL_{\ell y} ~\bs_\ell \\
  	   \sum_{\ell \geqslant 0} \tilde{\bd}_\ell\transpose ~\bL_{\ell z} ~\bs_\ell \\
  	 \end{bmatrix} \\
  \bH
  &= 
  \begin{bmatrix}
  	\sum_{\ell \geqslant 0} \frac{ \tilde{\bd}_\ell\transpose \{\bL_{\ell x},\bL_{\ell x}\} \bs_\ell}{2} &
  	\sum_{\ell \geqslant 0} \frac{ \tilde{\bd}_\ell\transpose \{\bL_{\ell x},\bL_{\ell y}\} \bs_\ell}{2} &
  	\sum_{\ell \geqslant 0} \frac{ \tilde{\bd}_\ell\transpose \{\bL_{\ell x},\bL_{\ell z}\} \bs_\ell}{2} \\
  	* &
  	\sum_{\ell \geqslant 0} \frac{ \tilde{\bd}_\ell\transpose \{\bL_{\ell y},\bL_{\ell y}\} \bs_\ell}{2} &
  	\sum_{\ell \geqslant 0} \frac{ \tilde{\bd}_\ell\transpose \{\bL_{\ell y},\bL_{\ell z}\} \bs_\ell}{2} \\
  	* &
  	* &
  	\sum_{\ell \geqslant 0} \frac{ \tilde{\bd}_\ell\transpose \{\bL_{\ell z},\bL_{\ell z}\} \bs_\ell}{2}
  \end{bmatrix} 
\end{align}
where we define $\tilde{\bd}_\ell = \bM_\ell\transpose( \bR_{\tilde{\ba}} )\bd_\ell$, $\{\bL_{\ell j}, \bL_{\ell k}\} \triangleq \bL_{\ell j} \bL_{\ell k} + \bL_{\ell k} \bL_{\ell j}$ is the anticommutator, and some terms of $\bH$ are omitted due to the symmetry of the Hessian. 
Eq.~\eqref{eqn:corr_quadraitc_approx} gives the quadratic approximation of the correlation function in closed form, which follows directly from the Lie algebra formulation of RSH-rot.
This example demonstrates an important advantage of the Lie algebra approach.

\providecommand{\sci}[2]{#1{\times}10^{#2}}

\providecommand{\sci}[2]{#1{\times}10^{#2}}

\section{Numerical Experiments}
\label{sec:experiments}

This section evaluates the proposed formulation numerically.
Three ways of computing the RSH-rot matrix $\bM_\ell$ derived in this paper are compared with the recursion of Ivanic and Ruedenberg~\cite{IvanicRuedenberg96,IvanicRuedenberg98}, which is used as the reference.
Four experiments are performed:
\begin{enumerate}
  \item the computation time per rotation;
  \item the accuracy of the matrices;
  \item rotations about a fixed axis, which exploit the separation of axis and angle in the anti-diagonal factorization;
  \item the closed-form derivatives of Section~\ref{sec:differentiation}.
\end{enumerate}

\subsection{Compared methods and protocol}
\label{sec:exp-methods}

\paragraph{Methods.}
\emph{Ivanic} starts from $\bM_1=\mathbf{P}\bR\mathbf{P}^\top$, with the permutation $\mathbf{P}$ that maps $(x,y,z)$ to $(y,z,x)$, and obtains each band from the previous one, so band $\ell$ requires all bands $1,\dots,\ell$.
\emph{Lie-exp} evaluates $\bM_\ell=\exp(\mathbf{A}_\ell)$ with $\mathbf{A}_\ell=\theta_a(\hat{a}_x\bL_{\ell x}+\hat{a}_y\bL_{\ell y}+\hat{a}_z\bL_{\ell z})$ by the general-purpose matrix exponential of MATLAB.
\emph{Lie-eig} evaluates the same exponential through the numerical eigendecomposition of $\mathbf{A}_\ell$, $\bM_\ell=\operatorname{Re}(\mathbf{V}\operatorname{diag}(e^{im\theta_a})\mathbf{V}^{-1})$.
\emph{Lie-AD}, the proposed method, uses the anti-diagonal factorization $\bM_\ell=\bF_\ell\,\bM_{\ell z}(\theta_a)\,\bF_\ell^\top$ of Section~\ref{sec:rsh_rot_eigenvectors}, with $\bF_\ell$ obtained from the ladder recurrence of Proposition~\ref{thm:rsh_eigenvector_recursive}.
In floating point the recurrence amplifies round-off by the ratios $\gamma_{\ell,k}/\gamma_{\ell,m}$ of the ladder coefficients, and $\lVert\bF_\ell^\top\bF_\ell-\bI\rVert_F$ grows to about $10^{-8}$ at $\ell=64$.
Each new pair of vectors is therefore re-orthogonalized against the previous ones (classical Gram--Schmidt, two passes), which keeps $\bF_\ell$ orthogonal to round-off level, below $\sci{8.6}{-15}$ for all tested bands, at a cost of $\mathcal{O}(\ell^3)$ per band.
The generators depend only on $\ell$; they are computed once and excluded from the timings.

\paragraph{Protocol.}
All experiments were run in MATLAB R2026a on Intel Core i9-11900F@2.50 GHz CPU with 32 GB RAM.
Rotations are drawn from the Haar measure on $SO(3)$ with a fixed seed, every method receives the same $3\times3$ matrices $\bR$, and $\ell\in\{1,2,3,4,6,8,12,16,24,32,48,64\}$.
Each method was first checked independently against the defining relation $\bY_\ell(\bR\bphat)=\bM_\ell(\bR)\bY_\ell(\bphat)$ at random points.
All four methods satisfy it to within $\sci{4.2}{-15}$ in the same basis, while any other choice of phase or sign convention gives residuals of order one.
The differences reported below are therefore numerical.

\subsection{Experiment 1: computation time}
\label{sec:exp-time}

\begin{figure}[ht]
  \centering
  \includegraphics[width=0.62\linewidth]{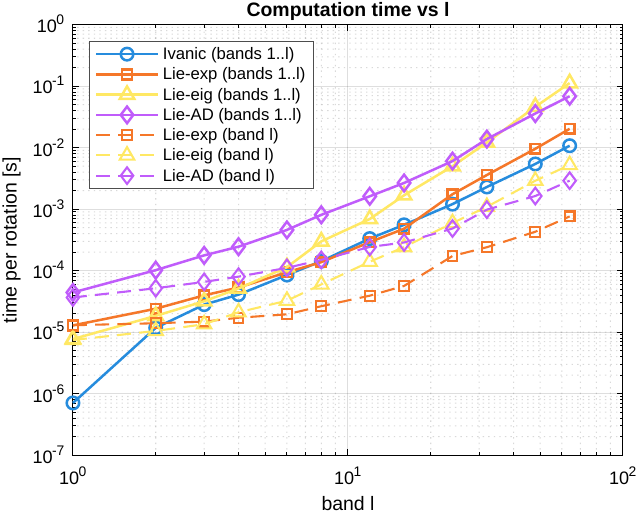}
  \caption{Computation time per rotation as a function of the band $\ell$. Solid lines: all bands $1,\dots,\ell$. Dashed lines: band $\ell$ only.}
  \label{fig:time}
\end{figure}

Each method processes a loop of random rotations timed with \texttt{timeit}.
Two settings are measured: all bands $1,\dots,\ell$, as required to rotate a complete expansion, and band $\ell$ only, which the recursion cannot compute without the lower bands (Fig.~\ref{fig:time}, Table~\ref{tab:timing}).

When all bands are required, the recursion is the fastest method.
Lie-exp is within a factor of two of it, and slightly faster for $8\leqslant\ell\leqslant16$.
Lie-AD is $4.8$ to $6.6$ times slower ($68.5$~ms against $10.8$~ms at $\ell=64$).
When a single band is required, the Lie-algebraic methods avoid the lower bands.
For $\ell=64$, Lie-exp and Lie-AD compute band $\ell$ respectively $14$ and $3.7$ times faster than the recursion, and Lie-AD is faster than the recursion from $\ell=12$ and faster than Lie-eig from $\ell=24$.

\begin{table}[ht]
  \centering
  \caption{Time per rotation in milliseconds, for all bands $1,\dots,\ell$ (left) and for band $\ell$ only (right).}
  \label{tab:timing}
  \begin{tabular}{r rrrr rrr}
    \toprule
     & \multicolumn{4}{c}{bands $1,\dots,\ell$} & \multicolumn{3}{c}{band $\ell$ only} \\
    \cmidrule(lr){2-5}\cmidrule(lr){6-8}
    $\ell$ & Ivanic & Lie-exp & Lie-eig & Lie-AD & Lie-exp & Lie-eig & Lie-AD \\
    \midrule
     4 &  0.04 &  0.05 &   0.05 &  0.25 & 0.02 & 0.02 & 0.08 \\
    16 &  0.56 &  0.48 &   1.68 &  2.67 & 0.06 & 0.24 & 0.29 \\
    32 &  2.29 &  3.59 &  12.10 & 13.80 & 0.24 & 1.09 & 0.98 \\
    64 & 10.80 & 20.20 & 112.00 & 68.50 & 0.76 & 5.22 & 2.90 \\
    \bottomrule
  \end{tabular}
\end{table}

\subsection{Experiment 2: accuracy}
\label{sec:exp-acc}

For $50$ random rotations per band, the discrepancy of each Lie-algebraic matrix from the recursion is recorded, together with the orthogonality residual $\lVert\bM_\ell^\top\bM_\ell-\bI\rVert_F$ of every method.
The residual requires no reference, since every exact RSH-rot matrix is orthogonal.

\begin{figure}[ht]
  \centering
  \includegraphics[width=\linewidth]{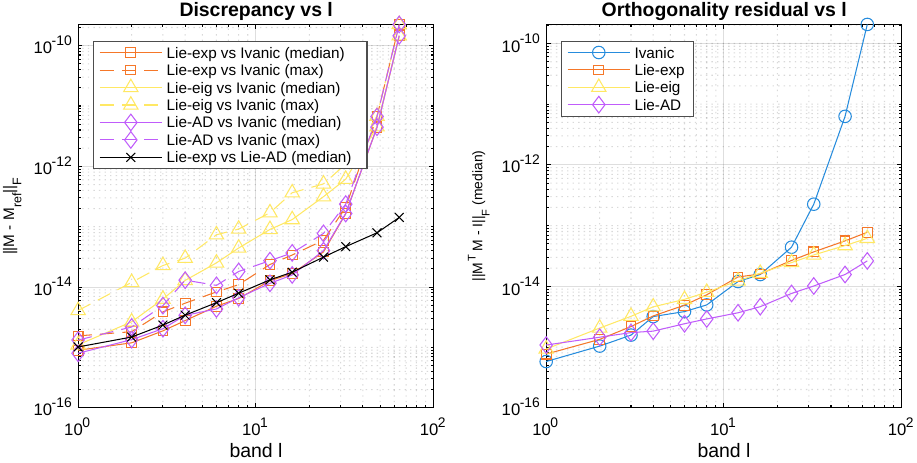}
  \caption{Left: discrepancy of the Lie-algebraic matrices from the Ivanic recursion (median and maximum over $50$ rotations) and between Lie-exp and Lie-AD. Right: median orthogonality residual of each method.}
  \label{fig:accuracy}
\end{figure}

Up to $\ell=16$ all methods agree with the recursion to within $\sci{1.3}{-13}$ (Fig.~\ref{fig:accuracy}).
At higher bands the three Lie-algebraic matrices differ from the recursion by the same amount, $\sci{1.5}{-10}$ at $\ell=64$, while Lie-exp and Lie-AD agree with each other to $\sci{1.4}{-13}$.
The orthogonality residuals identify the source.
The residual of the recursion grows to $\sci{2.1}{-10}$ at $\ell=64$, as its errors accumulate from one band to the next, whereas all Lie-algebraic methods, which compute each band independently, remain below $10^{-13}$.
Lie-AD produces the most nearly orthogonal matrices for every $\ell\geqslant4$, with a residual of $\sci{2.6}{-14}$ at $\ell=64$.

\subsection{Experiment 3: rotations about a fixed axis}
\label{sec:exp-axis}

In the anti-diagonal factorization, $\bF_\ell$ depends only on the axis, and the angle enters only through the sparse and explicit $\bM_{\ell z}(\theta)$.
Once $\bF_\ell$ is computed, rotations by any number of angles about the same axis require no further matrix function.
This is the case, for instance, when the heading of a sensor with known vertical direction must be estimated.
In particular, the correlation of eq.~\eqref{eqn:rsh_correlation} restricted to the axis becomes the trigonometric polynomial
\begin{equation}
  C(\theta)=\sum_{\ell=0}^{L}\mathbf{d}_\ell^\top\bM_\ell(\theta)\mathbf{s}_\ell
  = a_0+\sum_{m=1}^{L}\big(A_m\cos m\theta+B_m\sin m\theta\big),
  \label{eq:trig-poly}
\end{equation}
whose coefficients follow from $\bF_\ell^\top\mathbf{d}_\ell$ and $\bF_\ell^\top\mathbf{s}_\ell$.

The rotated coefficients $\bM_\ell(\theta_k)\mathbf{s}_\ell$ and the correlation $C(\theta_k)$ are computed for $K=360$ angles about a random axis, for bands $0,\dots,L$ of a three-lobe vMF mixture (eq.~\eqref{eqn:vmf_rsh_coeffs}), with closed-form reference values.
Lie-eig, which can also reuse its eigendecomposition for a fixed axis, is included as the natural competitor. We are generating kernels for $n = 3$, with variable $\kappa = [20, 80, 300]$ (ref. to eq. \eqref{eq:vmf_mixture_rsh}).
Times include the precomputation for the axis.

\begin{figure}[ht]
  \centering
  \includegraphics[width=\linewidth]{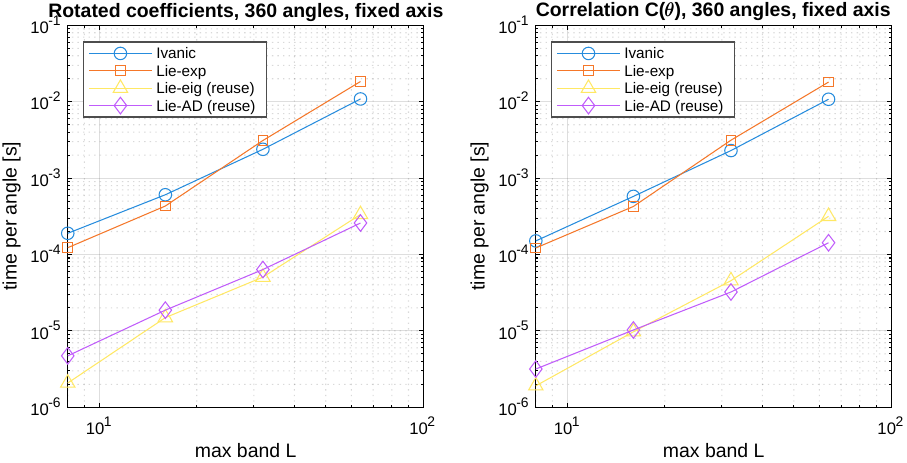}
  \caption{Rotations about a fixed axis: time per angle over a sweep of $K=360$ angles, including the precomputation for the axis. Left: rotated coefficients. Right: correlation $C(\theta)$.}
  \label{fig:fixed-axis}
\end{figure}

Reusing the axis reduces the cost by more than an order of magnitude (Fig.~\ref{fig:fixed-axis}).
For $L=64$, Lie-AD evaluates the correlation at the $360$ angles in $51.3$~ms, $76$ times faster than the recursion and $127$ times faster than Lie-exp.
Of this time, $49.0$~ms are spent computing $\bF_\ell$, so the whole sweep costs less than a single rotation of all bands by Lie-AD and about as much as five by the recursion.
Lie-AD is also faster than Lie-eig with reuse for $L=64$, by $2.2$ times for the correlation and $1.3$ times for the rotated coefficients.
At the same time it is up to seven times more accurate, since its factor $\bF_\ell$ is real and orthogonal.
Finally, the maximum of eq.~\eqref{eq:trig-poly} is found exactly by Newton iterations on the polynomial: the angle is recovered to working precision, against errors of $10^{-6}$--$10^{-5}$~rad for a grid search with parabolic refinement.

\subsection{Experiment 4: closed-form derivatives}
\label{sec:exp-deriv}

\begin{figure}[ht]
  \centering
  \includegraphics[width=\linewidth]{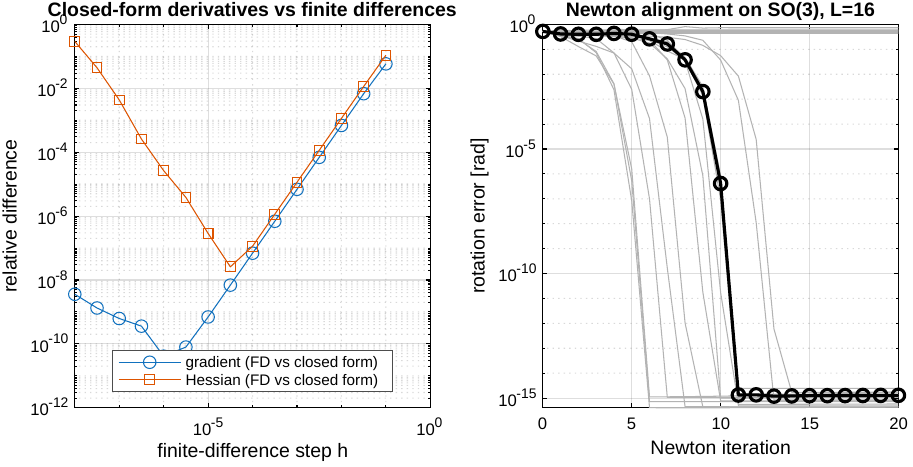}
  \caption{Left: relative difference between finite-difference and closed-form gradient and Hessian as a function of the step $h$. Right: rotation error of Newton alignment from $30^\circ$; individual trials in grey, median in black.}
  \label{fig:derivatives}
\end{figure}

The gradient $\bb$ and Hessian $\bH$ of the correlation function (Section~\ref{sec:differentiation}) are computed from the generators for $L=16$, at a rotation $20^\circ$ from the optimum.
They are compared with central finite differences of $C(\bR\exp(\bOmega_{\bomega}))$, which require $19$ evaluations of $C$.
The finite differences agree with the closed form at best to $\sci{4.1}{-11}$ for the gradient and $\sci{2.6}{-8}$ for the Hessian.
Their error grows for both smaller and larger steps (Fig.~\ref{fig:derivatives}, left), while the closed-form evaluation is exact up to round-off and $17$ times faster.

As an illustration, Newton iterations $\bomega=-\bH^{-1}\bb$, $\bR\leftarrow\bR\exp(\bOmega_{\bomega})$, with a backtracking line search, were started from $20$ rotations at $30^\circ$ from the optimum.
In $14$ trials they reached the optimum, with quadratic convergence in the final iterations down to a rotation error of $10^{-15}$~rad (Fig.~\ref{fig:derivatives}, right).
In the remaining trials they did not reach the optimum within $20$ iterations.
The remaining 6 trials converged to other local maxima of the correlation, so, as for any local method, a good initialization is required.

\subsection{Discussion}
\label{sec:exp-discussion}

For computing all bands of a single rotation, the Ivanic recursion remains the fastest method.
At high bands, however, it is the least accurate, while all Lie-algebraic methods stay at round-off level, and Lie-AD produces the most nearly orthogonal matrices.
The formulation is most advantageous when the problem has more structure than a single rotation.
It computes a single band without the lower ones, and it evaluates many rotations about the same axis at about the cost of one, with an exactly maximizable correlation.
It also gives exact gradients and Hessians at a fraction of the cost of finite differences.
The main limitation of Lie-AD is the stability of the ladder recurrence, which requires re-orthogonalization; exploiting the structure of $\bF_\ell$ to reduce its cost is left for future work.

\section{Conclusion}
\label{sec:conclusion}

We have presented a Lie-algebraic formulation of the rotation of real spherical harmonic coefficients.
We have derived in closed form the generators of the representation of $\mathfrak{so}(3)$ on each RSH band, from which the RSH-rot matrix follows through the exponential map.
As an alternative to a general-purpose matrix exponential, we have introduced an anti-diagonal factorization that confines the dependence on the rotation angle to an elementary rotation about the $z$ axis, the remaining factor being an orthogonal similarity matrix that depends on the axis alone.
The similarity matrix has been characterized by two independent routes, through the adjoint action on the Lie algebra and through the eigenvectors of the RSH-rot matrix; the latter are computed in real arithmetic by a ladder recurrence, which we have shown to require re-orthogonalization to be numerically stable.

The experiments show that the recursion of Ivanic and Ruedenberg remains the fastest method when all bands of a single rotation are required, but loses accuracy at high bands, whereas the Lie-algebraic methods remain at round-off level, with the anti-diagonal factorization producing the most nearly orthogonal matrices.
The Lie-algebraic methods are faster when a single band is required.
Since the similarity matrix depends only on the axis, many rotations about the same axis cost about as much as one, and the correlation of two RSH expansions along that axis becomes a trigonometric polynomial that can be maximized exactly.
Finally, the generators give exact gradients and Hessians, so rotation-dependent objectives of RSH series are differentiable at a fraction of the cost of finite differences.
Reducing the cost of the stable computation of the similarity matrix, by exploiting its structure, is left for future work.

\appendix
\section{Complex Spherical Harmonics and Conversion formulas}
\label{sec:appendix_csh_conversion}

The Complex Spherical Harmonic (CSH) of band $\ell$ and index $m$  with normalization is defined as 
\begin{equation}
  Y_{l}^{m}(\theta, \varphi) 
  {=}(-1)^m \sqrt{\frac{2\ell{+}1}{4 \pi}\frac{(\ell{-}m)!}{(\ell{+}m)!}}
     P_\ell^{m}(\cos\theta) ~\EXPSHORT{i m \varphi} 
  \label{eqn:csh_definition}
\end{equation}
for $m\geqslant0$, with $Y_\ell^{-m}=(-1)^m\,\overline{Y_\ell^{m}}$. The above equation incorporates the customary Condon--Shortley phase $(-1)^m$ explicitly, and $P_\ell^{m}(\cos\theta)$ is the associated Legendre function without that phase, as in eq.~\eqref{eqn:rsh_definition}. 
The argument expressed by colatitude $\theta$ and longitude $\varphi$ angles is also equivalently provided by the unit vector on unit sphere $\bphat \in \mathbb{S}^2$. 
Following a standard convention the CSH are denoted by superscript $m$, $Y_{\ell}^m$ whereas the RSH are addressed by subscript $m$, $Y_{\ell, m}$. 

The conversion formulas state how to write RSH with respect to CSH and vice versa. 
The algebraic expression is rather different in the cases of positive and negative $m$ and to avoid confusion it is convenient to assume that $0 < m \leqslant \ell$ and explicitly write the negative index as $-m$ in the functions $Y_{l}^{m}$ or $Y_{l, -m}$.
The formulas of RSH with respect to CSH are the following
\begin{align}
  Y_{\ell, -m}
  &= \frac{i}{\sqrt{2}} ~Y_{\ell}^{-m} - \frac{i ~(-1)^m}{\sqrt{2}} ~Y_{\ell}^{m} \\
  Y_{\ell 0} 
  &= Y_{\ell}^{0} \\
  Y_{\ell,  m}
  &= \frac{1}{\sqrt{2}} ~Y_{\ell}^{-m} + \frac{(-1)^m}{\sqrt{2}} ~Y_{\ell}^{m} 
\end{align}
The formulas of CSH with respect to RSH are
\begin{align}
	Y_{\ell}^{-m}
	&= - \frac{i}{\sqrt{2}} ~Y_{\ell, -m} + \frac{1}{\sqrt{2}} ~Y_{\ell, m}  \\
	Y_{\ell}^{0} 
	&= Y_{\ell 0} \\
	Y_{\ell}^{m}
	&= i ~\frac{(-1)^m}{\sqrt{2}} ~Y_{\ell, -m} + \frac{(-1)^m}{\sqrt{2}} ~Y_{\ell, m} 
\end{align}

\section*{Declarations}

\paragraph{Funding}
The authors did not receive support from any organization for the submitted work.

\paragraph{Competing interests}
The authors have no competing interests to declare that are relevant to the content of this article.

\paragraph{Data availability}
The datasets analyzed during the current study are available from the corresponding author on reasonable request. 

\paragraph{Code availability}
The MATLAB implementation of the methods and the scripts that reproduce all experiments and figures are available at \url{https://www.github.com/dlr1516/rimsh}.

\bibliographystyle{spmpsci}      %

\bibliography{roboticbiblio,biblio}

\end{document}